\documentclass[12pt]{article}
\usepackage{amssymb}
\usepackage{amsfonts}
\usepackage{amsmath}
\usepackage[pdftex]{graphicx}

\usepackage{amssymb,amsmath,amsfonts,bbm}
\usepackage[pdftex]{graphicx}
\usepackage{mathrsfs}
\usepackage{calligra}
\usepackage{multirow}
\numberwithin{equation}{section}

\newcommand{\PMB}[1]{
\leavevmode
\setbox0=\hbox{#1}%
\kern-0.02em\copy0\kern-\wd0
\kern0.02em\copy0\kern-\wd0
\kern-0.025em\raise0.0167em\box0
\kern-0.025em\raise0.0433em\box0}

\begin{document}

\begin{center}
{\bf \Large {Intermediate efficiency of some weighted goodness-of-fit statistics}}

\vspace{0.5cm}

\vspace{0.5cm}

{\bf Bogdan \'Cmiel}\\
{\it Faculty of Applied Mathematics, AGH University of Science and Technology,\\ Al. Mickiewicza 30, 30-059 Cracov, Poland}\\

{\it e-mail}: {\tt cmielbog@gmail.com}\\

\vspace{0.5cm}
{\bf Tadeusz Inglot}\\
{\it Faculty of Pure and Applied Mathematics, Wroc{\l}aw University of Science and Technology,\\
Wybrze\.ze Wyspia\'nskiego 27, 50-370 Wroc{\l}aw, Poland}\\
{\it e-mail}: {\tt Tadeusz.Inglot@pwr.edu.pl}\\

\vspace{0.5cm}
{\bf Teresa Ledwina}\\
{\it Institute of Mathematics, Polish Academy of Sciences,\\
ul. Kopernika 18, 51-617 Wroc{\l}aw, Poland}\\
{\it e-mail}: {\tt ledwina@impan.pl}\\

\end{center}
\vspace{0.5cm}

{\it Abstract}: 
This paper compares the Anderson-Darling and some Eicker-Jaeschke statistics to the classical unweighted Kolmogorov-Smirnov statistic. The goal is to provide a quantitative comparison of such tests and to study real possibilities of using them to detect departures from the hypothesized distribution that occur in the tails. 
This contribution covers the case when under the alternative a moderately large portion of probability mass is allocated towards the tails. It is demonstrated that the approach allows for tractable, analytic comparison between the given test and the benchmark, and for reliable quantitative evaluation of weighted statistics. Finite sample results illustrate the proposed approach and confirm the theoretical findings.
In the course of the investigation we also prove that a slight and natural modification of the solution proposed by Borovkov and Sycheva (1968) leads to a statistic which is a member of Eicker-Jaeschke class and can be considered an attractive competitor of the very popular supremum-type Anderson-Darling statistic.\\ 



\newpage

\noindent
{\bf \large {1. Introduction}}\\ 

\noindent
Weighted Kolmogorov-Smirnov-type goodness-of-fit tests have received a renewed interest in recent years; cf. Jager and Wellner (2004, 2007), Chicheportiche and Bouchaud (2012), Greenshtein and Park (2012), Charmpi and Ycart (2015), Gontscharuk et al. (2016), and Stepanova and Pavlenko (2018) for some illustration. 
A renaissance in research has to a large extent been driven by an application of a supremum version of the Anderson-Darling statistic in detecting sparse heterogenous mixtures, invented and developed by Donoho and Jin (2004, 2015). 
Obviously, weighted statistics of supremum-type are useful in many other problems as well. The renewed interest raises many unsolved questions for such structures; cf. the list of open problems on p. 2032 in Jager and Wellner (2007), and Section 5 in Ditzhaus (2018), for example. One of the questions concerns the power behavior of the considered statistics under nearby alternatives. Another one involves better understanding of the advantages and limitations of popular classes of nonparametric statistics, reconsidered recently in the context of detection of some mixtures. The aim of the present paper is to provide some tools and at least partial answers to these challenging questions.

For an exemplification of our approach, we study some selected Eicker-Jaeschke-type statistics and  compare them with the classical Kolmogorov-Smirnov and the integral Anderson-Darling statistics. We focus on uniformity testing and restrict our attention to two representatives of the class:
$$
{\cal S}_n = \sqrt n \sup_{0 < t < 1} \frac{|\hat F_n(t) - t|}{\sqrt{t(1-t)}}
\eqno(1.1)
$$
and its truncated variant
$$
{\cal E}_n = {\cal E}_n(\kappa_n)=\sqrt n \sup_{\kappa_n \leq t \leq 1-\kappa_n} \frac{|\hat F_n(t) - t|}{\sqrt{t(1-t)}},\;\;\;\kappa_n \in (0,1/2),\;\kappa_n\to 0\;\;\mbox{as}\;\;n \to \infty,
\eqno(1.2)
$$
where $\hat F_n(t)$ is the empirical distribution function of $n$ independent random variables with values in (0,1). ${\cal S}_n$ was proposed by Anderson and Darling (1952) while ${\cal E}_n$ is a consistent variant of the statistic 
$$
{\cal G}_n = {\cal G}_n(\kappa)=\sqrt n \sup_{\kappa \leq t \leq 1-\kappa} \frac{|\hat F_n(t) - t|}{\sqrt{t(1-t)}},\;\;\;\kappa \in (0,1/2),
\eqno(1.3)
$$
introduced and studied by Borovkov and Sycheva (1968).

Borovkov and Sycheva (1968) have shown that if the type I error tends to 0 slower than exponentially, as $n \to \infty$, then the uniform weight function $1/{\sqrt{t(1-t)}}$ ensures that ${\cal G}_n$ is asymptotically uniformly most powerful, in a certain sense, in some class of weighted statistics. A similar result for an exponentially decreasing type I error is contained in Borovkov and Sycheva (1970). Eicker (1979) and Jaeschke (1979) have obtained Darling-Erd\"{o}s-type results for 
${\cal S}_n$ and ${\cal E}_n$, under the null model, and suggested that ${\cal S}_n$ is sensitive in detecting moderate tails while, in contrast, the classical unweighted Kolmogorov-Smirnov test, say ${\cal K}_n$, is asymptotically sensitive in detecting changes in the central range of the null distribution. R\'ev\'esz (1982) provided some illustrative results supporting such statements, while Mason and Schuenemeyer (1983, 1992) defined and studied some formalization of the ability to detect central and local tail departures. 
They also studied a class of R\'enyi-type tests, being also weighted statistics, but with heavier weights than the uniform one. 
Jager and Wellner (2007) studied, among others, the optimal detection boundary of ${\cal S}_n$ for a sparse heterogenous mixture model. Ditzhaus (2018) extended the results in Jager and Wellner (2007) in many directions. 
Based on the findings of the two above mentioned papers, one sees that from the point of view of complete detectability of specific signals, a very large class of tests was shown to achieve the same completely detectable region, under very general signal models, as the very popular higher criticism test, related to the supremum-type Anderson-Darling statistic. It should be also strongly emphasized that all the above mentioned results on different forms of detectability were phrased in terms of the presence or absence of a power consistency under some convergent sequences of alternatives.

We would like to propose some quantitative results to study local power of some representatives of currently popular statistics from another perspective.
Namely, an interesting question is how many observations are needed for these tests to attain a given power lying in the interval $(0,1)$. Therefore, we shall compare the related numbers of observations via an appropriate asymptotic relative efficiency (ARE) notion. 
Moreover, we would like to show that careful introduction of the uniform weight, in a way proposed in (1.2), results in a stable and highly efficient solution. 
Surprisingly enough, this member of the Eicker-Jaeschke class has thus far received much less attention than ${\cal S}_n$. To complete the picture of the sup-type Anderson-Darling statistic, we also consider its integral variant.

Our approach to computing the efficiency of the considered statistics relies on a pathwise variant of Kallenberg's intermediate ARE.
The variant, elaborated in Inglot et al. (2018), is flexible enough to be applicable to some cases which lack high regularity. Weighted goodness-of-fit statistics, based on the classical empirical process, fall into this category. 
The characteristic features of the intermediate efficiency are: type I error tending to 0 slower than exponentially; local alternatives converging to the null distribution slower than $1/\sqrt n$; and, in contrast to the above mentioned developments on different forms of distinguishability, and non-degenerate asymptotic powers under local alternatives. The efficiency shares the advantages of Bahadur's and Pitman's approaches, but is much more widely applicable. 
In particular, the intermediate efficiency exploits moderate deviations of test statistic under the null model, instead of large deviations inherent in Bahadur's theory. For many weighted statistics large deviations are degenerate while moderate deviations are not.
For a more detailed discussion, see Inglot et al. (2018).

In our efficiency calculations the classical unweighted Kolmogorov-Smirnov statistic ${\cal K}_n$ shall play the role of a benchmark with respect to which other statistics shall be compared. Basically, to get the efficiency, one has to guarantee non-degenerate asymptotic powers of test statistics under a given sequence of alternatives and non-degenerate moderate deviations under the corresponding null model. The last question calls for example for using ${\cal M}_n = \sqrt {\log({\cal S}_n +1)}$ in place of ${\cal S}_n$. Sequences of alternatives are described in Section 2. In principle, they are defined via a fixed alternative distribution and a sequence of real parameters shrinking it to the null distribution. The efficiency allows for tractable analytic comparisons between two tests. 

We give a sufficient condition on tails of local sequences under which the intermediate efficiency of ${\cal E}_n$ with respect to ${\cal K}_n$ exists and is positive. 
Under this condition, ${\cal E}_n$ is always at least as efficient as 
${\cal K}_n$ and the efficiency of ${\cal E}_n$ with respect to ${\cal K}_n$ is always greater or equal to the efficiency ${\cal G}_n$ with respect to ${\cal K}_n$.
Moreover, we provide a sufficient condition, slightly stronger than that needed for ${\cal E}_n$, under which the efficiency of ${\cal M}_n$ with respect to ${\cal K}_n$ exists and is 0. In such situations, ${\cal E}_n$ does much better than ${\cal M}_n$, as a rule. Both sufficient conditions define local alternatives which do not shift too much mass towards one or two ends of $(0,1)$ and provide clear hints on which departures from the null model can or can not be detected by ${\cal E}_n$ and ${\cal M}_n$, respectively. Besides, the values of the efficiency nicely reflect the finite sample powers. We illustrate this in Section 8, where testing for the standard Gaussian distribution is considered. In Section 9 we study the case when the tails of the alternative are more heavy than they were assumed in Section 8. We compare there the above mentioned tests via simulations and state the result saying that so-called weak variant of the intermediate efficiency of ${\cal E}_n$ with respect to ${\cal K}_n$ is infinite. The outcomes, along with the results of Section 8, show that ${\cal E}_n$ with a relatively small smoothing parameter $\kappa_n$ is a well balanced solution working nicely under different kinds of tails of alternatives.

The structure of the paper is as follows: In Section 3 we restate slightly generalized results of Inglot and Ledwina (2006) related to the Kolmogorov-Smirnov statistic ${\cal K}_n$. Sections 4 and 5 collect necessary technical results on ${\cal M}_n$ and ${\cal E}_n$. Section 6 presents respective results on the integral Anderson-Darling statistic ${\cal I}_n$. Section 7 gives analytical formulas for the Kallenberg efficiencies of ${\cal M}_n$, ${\cal E }_n$, and ${\cal I}_n$ with respect to ${\cal K}_n$, and discusses the results. Section 8 reports outcomes of some simulation experiments. Section 9 contains some preliminary study of efficiency of ${\cal E}_n$ with respect to ${\cal K}_n$ under heavy-tailed alternatives.
We close with Section 10 containing some discussion of our results. All proofs are collected in the Appendix. 
\\

\noindent
{\bf \large {2. Testing problem and sequences of alternatives}}\\

\noindent
Throughout we rely on the setup and results of Inglot et al. (2018). 
As typical in the one-sample case, we denote the sample size by $n$ instead of $N$, as it was done in the general problem considered ibidem.
Let $X_1,...,X_n$ be independent random variables with continuous distribution function $F$. 
Denote by $F_0$ the null distribution function. Consider testing 
$$\mathbb{H}_0 : F=F_0$$
against the unrestricted alternative
$$
\mathbb{H}_1 : F \neq F_0.
$$

To introduce a class of sequences of alternatives approaching to $F_0$, consider first a fixed alternative $F_1$, a parameter $\vartheta_n \in (0,1)$, the combination 
$(1-\vartheta_n)F_0 + \vartheta_n F_1$ and its transformation to (0,1) via $F_0$. This yields the following alternative to the uniform distribution on (0,1)
$$
F_n^*(t)=[(1-\vartheta_n)F_0 + \vartheta_n F_1] \circ F_0^{-1} = t + \vartheta_n[F_1 \circ F_0^{-1}(t) - t].
\eqno(2.1)
$$
The function $F_1 \circ F_0^{-1}$ is called the comparison distribution function or the ordinal dominance curve. If $F_1$ is absolutely continuous with respect to $F_0$ then the density, say $f^*$, of $F_1 \circ F_0^{-1}(t)$ with respect to the Lebesgue measure on $(0,1)$ exists. The density is labeled as the comparison density, the relative density or the grade density. In terms of densities, (2.1) reads as $f_n^*(t)=(1-\vartheta_n){\bf 1}_{(0,1)}(t)+\vartheta_n f^*(t)$, where ${\bf 1}_{(0,1)}(t)$ stands for the uniform density on $(0,1)$. See Handcock and Morris (1999), and Thas (2010) for details.

The above motivates us to consider the observations from $[0,1]$, $\mathbb{H}_0$ : $F(t)=t,\;t \in (0,1)$, and
nearly null distribution functions of the form 
$$
F_n(t) = t + \vartheta_n A(t),\;\;\;t \in (0,1), 
\eqno(2.2)
$$
where $A(t)$ is continuous, $A(0)=A(1)=0,\;A \not\equiv 0$, while $\vartheta_n \to 0$ as $n \to \infty$. In many standard situations the function $A$ is absolutely continuous with a derivative $a$, which is unbounded. For an illustration see Section 8. This is in sharp contrast to the situation we considered in the two-sample problem, treated in Inglot et al. (2018).

In what follows, by $P_{\vartheta_n}$ we denote the probability measure related to $F_n$ in (2.2) while $P_0$ stands for the uniform distribution on $(0,1)$.
Moreover, $P_{\vartheta_n}^n$ and $P_0^n$ denote $n$ fold products of $P_{\vartheta_n}$ and $P_0$, respectively.\\

\noindent
{\large { 3. {\bf The intermediate slope of the classical Kolmogorov-Smirnov statistic} ${\cal K}_n$}}\\ 

\noindent
We have 
$$
{\cal K}_n = \sqrt n \sup_{0 <t <1} |\hat F_n(t) - t|,
\eqno(3.1)
$$
where $\hat F_n$ is the empirical distribution function of the sample. The intermediate slope of ${\cal K}_n$, under (2.2) with $\vartheta_n \to 0$ in such a way that $\sqrt n \vartheta_n \to \infty$, can be deduced from Inglot and Ledwina (2006). However, it should be noted that in that paper the corresponding sequences of alternatives were defined via densities. This forced an unnecessary assumption that the related $A$ should be absolutely continuous. Moreover, for convenience, it was assumed that $a=A'$ is bounded. Under (2.2) no extra assumptions are needed. For completeness, we restate here the corresponding results. In particular, (3.4), below, follows immediately from the proof of Theorem 6.1 in Inglot and Ledwina (2006). 

Define 
$$||A||_{\infty}=\sup_{0 <t< 1}|A(t)|\;\;\;\mbox{and}\;\;\;b_{{\cal K}}(P_{\vartheta_n}^n)=\sqrt n \vartheta_n||A||_{\infty}.\eqno(3.2)$$

\vspace{2mm}
\noindent
{\bf Proposition 1.} {\it For any positive $\{w_n\}$, such that $w_n\to 0$ and $n w_n^2 \to \infty$, as $n \to \infty$, it holds that 
$$
-\lim_{n \to \infty} \frac{1}{nw_n^2} \log P_0^n ({\cal K}_n \geq \sqrt n w_n)=c_{{\cal K}}=2
\eqno(3.3)
$$
and
$$
\lim_{n \to \infty} P_{\vartheta_n}^n \Bigl(\Big|\frac{{\cal K}_n}{b_{{\cal K}}(P_{\vartheta_n}^n)} -1\Big| \leq \epsilon\Bigr)=1
\eqno(3.4)
$$
for every $\epsilon > 0$. Consequently, the intermediate slope of ${\cal K}_n$ is $\; c_{{\cal K}}[b_{{\cal K}}(P_{\vartheta_n}^n)]^2=2n\vartheta_n^2||A||_{\infty}^2$.}\\

Note that for ${\cal K}_n$ we have the moderate deviations (3.3) in the full range of $w_n$'s and (3.4) holds without any further assumptions on $A$. As said before, ${\cal K}_n$ shall play the role of a benchmark procedure in our comparisons. 

In the next section we list some weighted variants of ${\cal K}_n$, which we shall further study, and present their moderate deviations under the null model. 
It should be emphasized that, in contrast to the benchmark procedure, the competitors do not need to have non-zero moderate deviations in the full range of $w_n$'s. This is very useful, as we shall see that it is a natural and an unavoidable restriction in the case of some weighted statistics. 

To calculate the efficiencies of weighted statistics, with respect to ${\cal K}_n$, we need for them some results analogous to (3.3) and (3.4) and, additionally, we have to identify sequences of alternatives for which asymptotic powers of these competitors of ${\cal K}_n$ are non-degenerate. 
These questions are solved in Sections 4 and 5. To get such asymptotic results, we shall consider some subclasses of functions $A$ in (2.2). The requirements are not very restrictive and many commonly used models fulfill them. \\

\noindent
{\bf {\large 4. Some weighted variants of ${\cal K}_n$ and their moderate deviations under $\mathbb{H}_0$}}\\ 

\noindent
In addition to the statistics ${\cal S}_n$ and ${\cal E}_n$, which are central in our study, for the purpose of some discussion we consider two additional statistics: ${\cal G}_n$, defined in (1.3), and 
$$
{\cal C}_n ={\cal C}_n(\tau)=\sqrt n \sup_{0 < t < 1} \frac{|\hat F_n(t) - t|}{[{t(1-t)}]^\tau},\;\;\;\tau \in (0,1/2),
\eqno(4.1)
$$
extensively investigated in the probabilistic literature; see Shorack and Wellner (1986) for some evidence.

For any of the above weighted statistics, say ${\cal T}_n$, we study for which sequences $\{w_n\}$, such that $w_n \to 0$ and $nw_n^2 \to \infty$, the limit 
$$
-\lim_{n \to \infty} \frac{1}{nw_n^2} \log P_0^n ({\cal T}_n \geq \sqrt n w_n)=c_{{\cal T}}
$$
exists. The number $c_{{\cal T}}$ is called the index of moderate deviations. Depending on whether $c_{{\cal T}} > 0$ or $c_{{\cal T}} = 0$, we speak of non-degenerate or degenerate moderate deviations. 

Obviously, the simplest solution is ${\cal G}_n$. For this statistic, similarly as for ${\cal K}_n$, moderate deviations exist and are non-degenerate in the whole range of $w_n$'s; cf. Lemma 1, below. As (3.3), the result is obtained by matching the KMT strong approximations and an asymptotic behavior of corresponding suprema of a weighted Brownian bridge. The last question is well studied, see Sec. II of Adler (1990) for some basic results and Lifschits (1995), Sec. 14, for further developments. The proof is skipped, as it is very similar to that for ${\cal K}_n$; cf. Inglot and Ledwina (1990) for details on ${\cal K}_n$.

The statistic ${\cal E}_n$ can be seen to be a refined variant of ${\cal G}_n$. In this case the situation is much more complex. 
Namely, if $\kappa_n$ tends to 0 relatively slowly, then, using again the strong approximation technique, we get non-degenerate moderate deviations. However, if the rate of convergence of $\kappa_n$ is too fast, then the index of moderate deviations is 0 for large class of sequences $\{w_n\}$; see Lemma 2. 
An even more extreme situation occurs in the case of ${\cal C}_n$, for which the moderate deviations are non-degenerate only for a very restricted class of sequences $\{w_n\}$; cf. (ii) of Lemma 4. For ${\cal S}_n$ the index of moderate deviations is 0 for all allowable sequences $\{w_n\}$'s; see (i) of Lemma 3. 
In such circumstances, similarly as in the case of the Bahadur approach to an efficiency, one can search for a monotonic function (or a sequence of functions), which, after imposing on a given statistic, leads to tails commensurable with that of ${\cal K}_n$. Obviously, such a monotonic transformation gives an equivalent test. It turns out that in the case of ${\cal S}_n$ the transformation $ x \to \sqrt{\log(1+x)}$ does the job and 
$$
{\cal M}_n = \sqrt {\log ({\cal S}_n + 1)}
\eqno(4.2)
$$
exhibits a quantifiable moderate deviation behavior. The result is due to Mason (1985); cf. (ii) of Lemma 3, below. Similarly, the second statement in (i) of Lemma 4 is due to Mason (1985). 

\vspace{5mm}
\noindent
{\bf Lemma 1.} {\it For any $w_n \to 0$, and such that $n w_n^2 \to \infty$, it holds that }
$$
-\lim_{n \to \infty} \frac{1}{nw_n^2} \log P_0^n ({\cal G}_n \geq \sqrt n w_n)=c_{{\cal G}}=1/2.
\eqno(4.3)
$$

\vspace{5mm}
\noindent
{\bf Lemma 2.} \\
(i) {\it Assume that $n\kappa_n \to \infty$. Then for any $w_n \to 0$, and such that $w_n/ \sqrt {\kappa_n} \to \infty$ 
it holds}
$$
-\lim_{n \to \infty} \frac{1}{nw_n^2} \log P_0^n ({\cal E}_n \geq \sqrt n w_n)=0.
$$
(ii) {\it Suppose $\liminf_{n\to\infty}n\kappa_n/\log^2n>0$. Then for any $w_n \to 0$, and such that $w_n=o(\sqrt{\kappa_n})$ and $n w_n^2/\log\log n \to \infty$, it holds that }
$$
-\lim_{n \to \infty} \frac{1}{nw_n^2} \log P_0^n ({\cal E}_n \geq \sqrt n w_n)=c_{{\cal E}}=1/2.
\eqno(4.4)
$$

\vspace{5mm}
\noindent
{\bf Lemma 3.} \\
(i) {\it If $w_n \to 0$ and $nw_n^2 \to \infty$ then}
$$
-\lim_{n \to \infty} \frac{1}{nw_n^2} \log P_0^n ({\cal S}_n \geq \sqrt n w_n)=0.
$$
(ii) {\it For any $w_n\to 0$, and such that $n w_n^2/\log \log n \to \infty$, we have}
$$
-\lim_{n \to \infty} \frac{1}{nw_n^2} \log P_0^n ({\cal M}_n \geq \sqrt n w_n)=c_{{\cal M}}=2.
\eqno(4.5)
$$

\vspace{5mm}
\noindent 
{\bf Lemma 4.}\\
(i) {\it Suppose that $w_n \to 0$ and $nw_n^2/\log n \to \infty$. Then for any $\tau \in (0,1/2)$}
$$
-\lim_{n \to \infty} \frac{1}{nw_n^2} \log P_0^n ({\cal C}_n \geq \sqrt n w_n)=0\;\;\;{and}\;\;\;
-\lim_{n \to \infty} \frac{1}{nw_n^2} \log P_0^n (\sqrt{\log ({\cal C}_n + 1)} \geq \sqrt n w_n)=\frac{1}{\tau}.
$$
(ii) {\it For any $w_n \to 0$, and such that $nw_n^2 \to \infty,\; w_n =o(\sqrt{\log n /n})$, we have for $\tau \in (0,1/2)$}
$$
-\lim_{n \to \infty} \frac{1}{nw_n^2} \log P_0^n ({\cal C}_n \geq \sqrt n w_n)=c_{{\cal C}}=2^{1-4 \tau}.
\eqno(4.6)
$$
\\

\noindent
{\bf Remark 1.} With probability 1 it holds that
$$
{\cal E}_n=\sqrt{n}\max\Big\{\max_{i:F_0(X_{(i)})\in [\kappa_n,1-\kappa_n]}\frac{\max\{|F_0(X_{(i)})-\frac{i}{n}|,|F_0(X_{(i)})-\frac{i-1}{n}|\}}{\sqrt{F_0(X_{(i)})(1-F_0(X_{(i)}))}}, \frac{T_n}{\sqrt{\kappa_n(1-\kappa_n)}}\Big\},
$$
where $T_n=\max\{|I_{1n}/n-1/n-\kappa_n|,|I_{2n}/n-1+\kappa_n|\}$, $I_{1n}=\min\{1\leq i\leq n+1: F_0(X_{(i)})>\kappa_n\},\;I_{2n}=\max\{0\leq i\leq n:F_0(X_{(i)})<1-\kappa_n\}$, $X_{(1)} \leq ... \leq X_{(n)}$ are order statistics of the sample $X_1,..., X_n$ while for convenience we additionally set $F_0(X_{(0)})=0,\;F_0(X_{(n+1)})=1$. 
Lemma 2 (ii) and the above exhibit that abandoning some fraction 
of smallest and largest transformed observations in the sample allows for non-degenerate moderate deviations when using the uniform weight. The above shows also that the construction of the statistic ${\cal E}_n$ follows a similar idea as the modified higher criticism statistic $HC_n^{+}$ defined in Section 3 of Donoho and Jin (2004), where a slightly smaller fraction of smallest transformed observations was abandoned. 
Some simulated powers of $HC_n^{+}$ are reported and discussed in Li and Siegmund (2015).\\

The proof of Lemma 2 is provided in the Appendix. Also there, we justify the index $0$ appearing in Lemma 3 and 4. The statement (ii) of Lemma 4 is a consequence of Proposition 2.5 in Inglot and Ledwina (1993). As mentioned earlier, (4.5) and the moderate deviations for $\sqrt{\log ({\cal C}_n + 1)}$ follow from Mason (1985).
The above shows that even such standard weighted statistics behave very differently, and this illustrates the ``irregularities'', we mentioned in Section 1. Anyway, 
for each of the considered examples there are sequences $\{w_n\}$ for which the respective index of moderate deviations is positive. This makes it possible to apply the pathwise variant of intermediate efficiency elaborated in Inglot et al. (2018). The next step in this direction is to study the asymptotic behavior of the statistics under sequences of alternatives. This question is studied below. To avoid repetitions of similar statements, we restrict our attention to presenting in full form only the respective results on ${\cal E}_n $ and ${\cal S}_n $.\\

\noindent
{\bf {\large 5. An asymptotic behavior of ${\cal E}_n$ and ${\cal S}_n$ under sequences of alternatives and their intermediate slopes}}\\ 

\noindent
We follow the scheme and notation of the definition of the pathwise variant of intermediate efficiency elaborated in Inglot et al. (2018). Therefore, we consider a particular sequence $\{\theta_n\},\;\theta_n \in (0,1),$ where $\theta_n \to 0$, as $n\to \infty$, and the related $F_n$ in (2.2), is given by
$$
F_n(t)=t+\theta_n A(t),\;\;\;t \in (0,1),
\eqno(5.1)
$$
where $A(t)$ is continuous and $A(0)=A(1)=0,\;A\not\equiv 0$. As in Section 2, we set $P_{\theta_n}$ for the distribution of $F_n$ and $P_{\theta_n}^n$ for its $n$-fold product. Additionally, introduce
$$
A^*(t)=\frac{A(t)}{\sqrt{t(1-t)}}.
\eqno(5.2)
$$

In the case of ${\cal E}_n$, assume that $A$ satisfies
$$ 
\lim_{t\to 0^+} A^*(t)=\lim_{t\to 1^-}A^*(t)=0, 
\eqno (5.3)
$$
where $A^*$ is defined in (5.2). Then there exists $\delta=\delta_{\cal E}(A) \in (0,1/2)$ such that 
$$
\sup_{t\notin [\delta,1-\delta]}|A^*(t)|=\frac{1}{2}\sup_{0<t<1}|A^*(t)|.
\eqno(5.4)
$$
Set
$$
b_{\cal E}(P_{\theta_n}^n)=\sqrt n \theta_n \sup_{0 < t < 1}|A^*(t)|.
\eqno(5.5)
$$

Throughout $\Phi(w),\;w \in \mathbb{R}$, stands for the standard normal distribution function.\\

\vspace{2mm}
\noindent
{\bf Theorem 1}. {\it Consider (5.1) with $A(t)$ satisfying (5.3) and $\theta_n \in (0,1),\;\theta_n = o(\sqrt {\kappa_n})$, and $n \theta_n^2/\log\log n \to \infty$. Then\\}
(i)\quad $\displaystyle \limsup_{n\to\infty} P_{\theta_n}^n({\cal E}_n-b_{\cal E}(P_{\theta_n}^n)\leqslant w)\leqslant E_2(w),\;\;w\in \mathbb{R}$;\\
(ii)\quad $\displaystyle \liminf_{n\to\infty} P_{\theta_n}^n({\cal E}_n-b_{\cal E}(P_{\theta_n}^n)\leqslant w)\geqslant E_1(w),\;\;w>0,$\\
\noindent
{\it where $E_2(w)=\Phi(w)$ is the standard normal distribution function, $E_1(w)$ is the distribution function of $\;\sup_{[\delta,1-\delta]}\bigl\{|B(t)|/\sqrt{t(1-t)}\bigr\}$ with $\delta$ defined in (5.4), while $B$ is a Brownian bridge.}\\
\noindent
{\it Hence}, ${\cal E}_n/ b_{\cal E}(P_{\theta_n}^n) \stackrel {P_{\theta_n}^n}{\longrightarrow} 1 $,
{\it and the intermediate slope of ${\cal E}_n$ under $\{P_{\theta_n}\}$ has the form\linebreak[4]
$c_{\cal E} [b_{\cal E}(P_{\theta_n}^n)]^2,$ where $c_{\cal E}=1/2.$}\\

\vspace{2mm}
\noindent {\bf Remark 2.} In the case of ${\cal G}_n$ an analogue of Theorem 1 holds true for any $A$ in (2.2). The only difference is that in the description of $E_1(w)$ one should use $\kappa$ in the place of $\delta$. Hence we get the following: For (5.1) with $\theta_n \in (0,1),\;\theta_n \to 0$, and $n \theta_n^2 \to \infty$, the intermediate slope of ${\cal G}_n$, under $\{P_{\theta_n}\}$, has the form
$c_{\cal G} [b_{\cal G}(P_{\theta_n}^n)]^2,$ where $c_{\cal G}=1/2,$ while $b_{\cal G}(P_{\theta_n}^n)=\sqrt n \theta_n \sup_{ \kappa \leq t \leq 1-\kappa}|A^*(t)|$. A comparison of $b_{\cal G}(P_{\theta_n}^n)$ with $b_{\cal E}(P_{\theta_n}^n)$ supports the statement that ${\cal E}_n$ is a natural refinement of ${\cal G}_n$.\\

We have also considered an analogue of Theorem 1 for ${\cal C}_n$ with fixed $\tau \in (0,1/2)$. The result, together with Lemma 4 (ii), shows that the intermediate slope of ${\cal C}_n$ is smaller than the related slope of ${\cal E}_n$. Hence, under fixed $\tau$, ${\cal C}_n$ is less efficient than ${\cal E}_n$. Therefore, we skip the presentation of the relevant details. 

We have also derived the intermediate slope of a recent modification of ${\cal S}_n$ introduced by Stepanova and Pavlenko (2018). The results do not differ substantially from these on ${\cal S}_n$. Therefore, we present here our results only for the classical case of ${\cal S}_n$.\\

For ${\cal M}_n=\sqrt{\log({\cal S}_n+1)}$ suppose that $A(t)$ satisfies
$$
\sup_{t \in (0,1)}\frac{|A(t)|}{\;[t(1-t)]^{1-\varpi}} < \infty \;\;\;\mbox{for some}\;\;\;\varpi \in [0,1/2).
\eqno(5.6)
$$
The assumption (5.6) implies that there exists $\delta=\delta_{\cal M}(A) \in (0,1/2)$ such that
$$
\sup_{t \notin [\delta,1-\delta]}|A^*(t)| = \frac{1}{2}\sup_{0<t<1}|A^*(t)|,
\eqno(5.7)
$$ 
where $A^*(t)$ is defined in (5.2). 
In terms of an alternative $F_1(w)$ in (2.1), the condition (5.6) means that 
$$
\sup_{w \in \mathbb{R}} \frac{|F_1(w)-F_0(w)|}{{\bigl\{F_0(w)[1-F_0(w)]\bigr\}}^{1-\varpi}} < \infty \;\;\;\mbox{for some}\;\;\;\varpi \in [0,1/2).
$$
Put
$$b_{\cal M}(P_{\theta_n}^n) = \sqrt { \log\Bigl(\sqrt n \theta_n \sup_{0<t<1}|A^*(t)|\Bigr)}=\sqrt{\log b_{\cal E}(P_{\theta_n}^n)}.
\eqno(5.8)
$$

\vspace{3.5mm}
\noindent 
{\bf Theorem 2}. {\it Suppose that $A(t)$ satisfies (5.6) with some $\varpi \in [0,1/2)$. 
 Consider (5.1) with $\theta_n \in (0,1),\; \theta_n=o(n^{-\varpi}) $ and $(\log n \theta_n^2)/\log \log n \to \infty$ as $n \to \infty.$ Then }\\
(i)\quad $\displaystyle \limsup_{n\to\infty} P_{\theta_n}^n({\cal S}_n-b_{\cal E}(P_{\theta_n}^n)\leqslant w)\leqslant S_2(w),\;\;w\in \mathbb{R}$;\\
(ii)\quad $\displaystyle \liminf_{n\to\infty} P_{\theta_n}^n({\cal S}_n-b_{\cal E}(P_{\theta_n}^n)\leqslant w)\geqslant S_1(w),\;\;w>0,$\\
\noindent
{\it where $S_2(w)=\Phi(w)$, $S_1(w)$ is the distribution function of $\;\sup_{[\delta,1-\delta]}\bigl\{|B(t)|/\sqrt{t(1-t)}\bigr\}$ with $\delta$ defined in (5.7), $b_{\cal E}(P_{\theta_n}^n)$ is defined in (5.5), while $B$ is a Brownian bridge.}\\
\noindent
{\it Hence}, ${\cal M}_n/ b_{\cal M}(P_{\theta_n}^n) \stackrel {P_{\theta_n}^n}{\longrightarrow} 1 $,
{\it and the intermediate slope of ${\cal M}_n$ under $\{P_{\theta_n}\}$ has the form
$c_{\cal M} [b_{\cal M}(P_{\theta_n}^n)]^2,$ where $c_{\cal M}=2.$}\\

\vspace{2mm}
\noindent {\bf Remark 3.} The restriction (5.6) on $A$, imposed in Theorem 2, is obviously stronger than the related condition (5.3) needed for ${\cal E}_n$.
When $A$ is absolutely continuous with a derivative $a$ and for some $\epsilon\in [0,1/2)$ it holds that
$\displaystyle \limsup_{t\to 0^+}t^{\epsilon} |a(t)|< \infty$ and $\displaystyle \limsup_{t\to 1^-} (1-t)^{\epsilon}|a(t)|<\infty$ then the condition (5.6) is satisfied with $\varpi=\epsilon$. In particular, when $a$ is bounded then (5.6) holds with $\varpi =0$. The case $\varpi\in (0,1/2)$ admits unbounded $a$. 

Consider the alternative (5.1) with $A$ of the form $A(t)=t^{\delta}-t,\;\delta\in (0,1/2)$.
Then (5.3) and (5.6) do not hold. This $A$ corresponds to a heavy-tailed departure.
When the null distribution is $F_0(x)=\Phi(x)$ then such $A$ corresponds to the Lehmann (1953) alternative $F_1(x)= \Lambda(x;\delta)=[\Phi(x)]^{\delta}$ in (2.1). For further discussion of some examples see Sections 8 and 9.\\

\noindent
{\bf {\large 6. The integral Anderson-Darling statistic ${\cal I}_n$ and the related asymptotic results}}\\

\noindent
Set
$$
{\cal I}_n = \Bigl\{n \int_0^1 \frac{[\hat F_n(t) - t]^2}{t(1-t)} dt\Bigr\}^{1/2}.
\eqno(6.1)
$$
By Proposition 2.2 and Remarks 2.2-2.4 in Inglot and Ledwina (1993) we infer the following.\\

\noindent
{\bf Lemma 5.} {\it For any $w_n \to 0$, and such that $n w_n^2 \to \infty$, it holds that }
$$
-\lim_{n \to \infty} \frac{1}{nw_n^2} \log P_0^n ({\cal I}_n \geq \sqrt n w_n)=c_{{\cal I}}=1.
\eqno(6.2)
$$

\vspace{2mm}
Now, consider alternatives of the form (5.1), with $A$ such that for some $\ell \in (0,1/2)$ it holds that
$$ 
\int_0^1\frac{|A(t)|^{2\ell}}{t(1-t)}dt<\infty.
\eqno (6.3)
$$
Observe that under (6.3) for $A^*(t)$ defined in (5.2) it holds that
$$
||A^*||_2 =
\Bigl\{\int_0^1 \frac{A^2(t)}{t(1-t)} dt\Bigr\}^{1/2} < \infty.
$$

Note that both conditions (5.3) and (5.6) imply (6.3). Asymptotic behavior of ${\cal I}_n$ under the sequence of alternatives (5.1) with $A$ satisfying (6.3) is described below. \\

\noindent
{\bf Theorem 3.} {\it Suppose $A(t)$ satisfies (6.3). Consider ${P_{\theta_n}}$ obeying (5.1) with $\theta_n\in(0,1),$ and such that $\;\theta_n \to 0$, $ n \theta_n ^2 \to \infty$ as $n\to\infty$. Then
$$
\lim_{n \to \infty} P_{\theta_n}^n \bigl({\cal I}_n -\sqrt n \theta_n ||A^*||_2 \leq w\bigr) = \Phi\Bigl(\frac{w||A^*||_2}{\rho_A}\Bigr),\;\;w\in \mathbb{R},
\eqno(6.4)
$$
where}
$$
\rho^2_A= \int_0^1 \int_0^1 \frac{[\min\{s,t\} -st]A(s)A(t)}{s(1-s)t(1-t)} ds dt.
$$
{\it Hence, the intermediate slope of ${\cal I}_n$ has the form} 
$$
c_{\cal I}[b_{\cal I}(P_{\theta_n}^n)]^2\;\;\; \mbox{where}\;\;\;c_{\cal I}= 1,\;\;b_{\cal I}(P_{\theta_n}^n)= \sqrt n \theta_n ||A^*||_2.
\eqno(6.5)
$$

\vspace{2mm}
The result (6.4) was reported in Inglot et al. (2000) for the case of $A(t)$ absolutely continuous with a bounded derivative $a(t)=A'(t)$. Its proof was very briefly sketched in Inglot et al. (1998). 
Here, for completeness, we provide detailed justification of (6.4). In fact, a result like (6.4) with the corresponding (6.3) can be immediately generalized to Hilbertian norms on $D[0,1]$ imposed on the empirical process. We omit the details.
Such a result, along with the technique developed in Inglot and Ledwina (1993), allows us to calculate intermediate slopes of a family of integral test statistics. \\

\noindent
{\bf Remark 4}. Assume that $A$ in (2.2) is absolutely continuous and $a=A'$. If $a \in L_r(0,1)$ for some $r > 1$ then (6.3) holds. If $a \in L_2(0,1)$ then (5.3) is satisfied. In the case $a \in L_r(0,1)$ for some $r >2$ we have (5.6) with $\varpi \geq 1/r$.

For testing $F_0(x)=\Phi(x)$ consider the alternative distribution function $F_1$, parametrized by $\zeta >0$, and given by $F_1(x)=\Pi (x;\zeta)$, where $\Pi (x;\zeta)=|x|^{-\zeta}/2$ if $x < -1$, $\Pi (x;\zeta)=1/2$ if $-1 \leq x \leq 1$, and $\Pi (x;\zeta)=1-x^{-\zeta}/2$ if $x >1$. $F_1$ is a member of the symmetric Pareto family considered in Grabchak and Samorodnitsky (2010). Such an $F_1$, via (2.1), corresponds to $A(t)=\Pi(\Phi^{-1}(t);\zeta) - t$ in (2.2). A simple calculation shows that (6.3) is satisfied for $\ell >2/\zeta$ when $\zeta>4$ while (5.3) does not hold for any $\zeta > 0$. Moreover, for each $\zeta >0$ it holds that $a=A' \notin L_r(0,1)$ for any $r >1$. In such a sense, $F_1$ has the heaviest possible tails which can appear in (2.1) when $F_0(x)=\Phi(x)$.\\

\noindent
{\bf {\large 7. Intermediate efficiencies of $\;{\cal G}_n,\;{\cal E}_n,\; {\cal I}_n,\;$ and $\;{\cal M}_n\;$ with respect to $\;{\cal K}_n$}}\\

\noindent
Exploiting the results collected in Sections 2 - 6 and using Theorem 1 from Inglot et al. (2018), we immediately obtain the following results.\\

\noindent
{\bf Theorem 4.} {\it Consider a sequence of alternatives $\{P_{\theta_n}\}$ defined by (5.1) with $n\theta_n^2\to\infty$.\\

\noindent (i) The intermediate efficiency of ${\cal G}_n$ with respect to ${\cal K}_n$, under the sequence $\{P_{\theta_n}^n\}$, exists and equals
$$
e_{{\cal G}{\cal K}} = e_{{\cal G}{\cal K}}(\kappa)=\frac{\sup_{\kappa\leq t\leq 1-\kappa}[A^*(t)]^2}{4||A||_{\infty}^2}; \eqno (7.1)$$

\noindent (ii) Suppose $\liminf_{n\to\infty}n\kappa_n/\log^2n>0$. If $A$ satisfies (5.3) and $\theta_n=o(\sqrt{\kappa_n}),\;n\theta_n^2/\log\log n\to\infty$. Then the intermediate efficiency of ${\cal E}_n$ with respect to ${\cal K}_n$ under the sequence $\{P_{\theta_n}^n\}$, exists and equals
$$
e_{{\cal E}{\cal K}} = \frac{\sup_{0 < t < 1}[A^*(t)]^2}{4||A||_{\infty}^2}; \eqno (7.2)$$

\noindent (iii) If $A$ satisfies (6.3) then the intermediate efficiency of ${\cal I}_n$ with respect to ${\cal K}_n$ under the sequence $\{P_{\theta_n}^n\}$, exists and equals
$$
e_{{\cal I}{\cal K}} = \frac{||A^*||_2^2}{2||A||_{\infty}^2}. \eqno (7.3)$$}

\vspace{3mm}
\noindent
{\bf Theorem 5.} 
{\it Consider a sequence of alternatives $\{P_{\theta_n}\}$ defined by (5.1) with $A$ satisfying (5.6) for some $\varpi \in [0,1/2)$ and $\theta_n=o(n^{-\varpi}),\;
(\log n\theta_n^2)/\log\log n\to\infty$. \\
Then the intermediate efficiency of ${\cal M}_n$ with respect to ${\cal K}_n$, under the sequence $\{P_{\theta_n}^n\}$, exists and equals
$$
e_{{\cal M}{\cal K}} = \lim_{n \to \infty} \frac{c_{\cal M} [b_{{\cal M}}(P_{\theta_n}^n)]^2}{c_{\cal K} [b_{\cal K}(P_{\theta_n}^n)]^2} = 0.
\eqno(7.4)
$$}

\vspace{2mm}
\noindent 
{\bf Remark 5.} We have chosen ${\cal K}_n$ as a benchmark since, first of all, it seems to be a natural reference statistic when some weighting is considered. Moreover, in view of the approach elaborated in Inglot et al. (2018), it is applicable in such a role since it obeys moderate deviations in the full range. Alternatively, in view of Lemmas 1 and 5, ${\cal G}_n$ and ${\cal I}_n$ can be used as benchmarks, as well. 
Perhaps the most natural candidate for a benchmark procedure could be the Neyman-Pearson test statistic for uniformity against $F_n$, cf. (2.2), defined when $A$ is absolutely continuous with derivative $a$. To justify such a choice, again one should know that moderate deviations for this statistic hold 
for all sequences $\{w_n\}$ such that $w_n \to 0$ and $nw_n^2 \to \infty.$
 This is the case when $a$ is bounded. However, for unbounded $a$ such a question seems to remain open. Results of Merlev\`ede and Peligrad (2009) suggest that for unbounded $a$ the speed $a_n=1/nw_n^2$, using their and our notations, needs to be adjusted to $\vartheta_n$.
\\

\vspace{2mm}
\noindent 
{\bf Remark 6.} The results (7.1), (7.2) and (7.3) show that, under appropriate assumptions, the sample sizes needed for the Kolomogorov-Smirnov test to be, given $\{P_{\theta_n}\}$, as good as the tests based on ${\cal G}_n$, ${\cal E}_n$ and ${\cal I}_n$, respectively, are equal approximately to $n e_{{\cal G}{\cal K}}$ $n e_{{\cal E}{\cal K}}$ $n e_{{\cal I}{\cal K}}$, respectively. Thus, they are approximately proportional to $n$. 

The relation (7.4) reveals that, under (5.6), for the Kolmogorov-Smirnov test the sample size sufficient to attain, given $P_{\theta_n}$, the power as good as that of the test based on ${\cal M}_n$ is of smaller order than $n$. A similar result to (7.4) can be formulated on $e_{{\cal M}{\cal I}}$.

The statement (7.4) deserves some more detailed comments. First of all, it should be emphasized that the intermediate efficiency concerns the situation when asymptotic powers of the corresponding tests are kept in $(0,1)$. 
Therefore, the result (7.4) does not contradict consistency of ${\cal M}_n$ under fixed or convergent alternatives. 
Observe that our approach exhibits that the functions $b_{\cal M}(\cdot)$ and $b_{\cal K}(\cdot)$, defining the intermediate slopes, are related to the respective shifts in the limiting theorems, which ensure non-degenerate asymptotic powers. Since $b_{\cal M}(\cdot) \ll b_{\cal K}(\cdot)$, it can be expected that, in a finite sample comparison, the power function of ${\cal M}_n$ should be much smaller than the corresponding power function of ${\cal K}_n$. This tendency is quantitatively measured by the intermediate slopes and the intermediate efficiency.
Since in the intermediate approach the alternatives are not very close to the null one and the levels do not decrease very fast, we can expect that a similar tendency shall be seen in empirical powers under fixed alternatives, which satisfy (5.6). In Section 8.2 we present a small simulation study which confirms such intuitions. 

Next, compare (7.4) with, consistent with it, findings of Lockhart (1991). In that paper it was shown that, under usual types of contiguous alternatives, the power and the level of ${\cal S}_n$ have the same limit, and the related ARE of the test with respect to the corresponding Neyman-Pearson test (NP) 
is 0. The same conclusion holds true for the ARE of ${\cal S}_n$ with respect to any other test with a nonzero asymptotic efficiency relative to the NP test.
In our opinion, in this application, the intermediate approach, resulting in non-zero shift, explains better observed empirical powers of ${\cal S}_n$ than the conclusion on the shift 0 under Pitman's approach.

The results on ${\cal S}_n$ in Lockhart (1991) were formulated in the case when $A$ is absolutely continuous and the corresponding function $a=A'$ belongs to $L_2(0,1)$. This assumption is standard in the classical approach to investigation of an asymptotic power and the asymptotic relative efficiency of tests under alternatives of order $1/\sqrt n$. For an illustration see the insightful results on ${\cal K}_n$ proved by Milbrodt and Strasser (1990), and Janssen (1995).
On the other hand, note that we have shown that the intermediate slope of ${\cal K}_n$ is well defined for any $a \in L_1(0,1)$. 
 This opens a possibility of comparisons of some competitors to ${\cal K}_n$ for some interesting alternatives with $a \notin L_2(0,1)$. Typically, alternatives with heavy tails lead, via $F_1 \circ F_0^{-1}$, to a corresponding $a \notin L_2(0,1)$. Tails heavier than Gaussian are common in many current applications. For related discussion, see Cont (2001).
Examples of such alternatives along with some preliminary results on a weak variant of the intermediate efficiency are presented in Section 9. It turns out that in such a setting the asymptotic behavior of ${\cal E}_n$ changes dramatically. Namely, the weak intermediate efficiency of ${\cal E}_n$ with respect to ${\cal K}_n$ is infinite. In light of recent results on the intermediate efficiency of the Neyman-Pearson statistic ${\cal V}_n$ with respect to ${\cal K}_n$, in the case when $a \in L_p(0,1),\;p \in (1,2),$ contained in Inglot (2019), this is not surprising. It turns out that, in contrast to ${\cal E}_n$, ${\cal K}_n$ is completely inefficient in such situations.\\

\vspace{2mm}
\noindent 
{\bf Remark 7.} 
An easy calculation shows that $e_{{\cal E}{\cal K}}\geq 1$, $e_{{\cal E}{\cal K}}\geq e_{{\cal G}{\cal K}}$ and $e_{{\cal I}{\cal K}}\leq 2e_{{\cal E}{\cal K}}$ for any $A$ satisfying (5.3). Moreover, $e_{{\cal G}{\cal K}}$ can be arbitrarily close to 0 (take $A'(t)=a(t)={\bf 1}_{[0,(\kappa+\delta)/2]}(t)-{\bf 1}_{((\kappa+\delta)/2,\kappa+\delta]}(t)$ for small $\delta>0$ , where ${\bf 1}_E$ denotes the indicator of the set $E$). On the other hand, $e_{{\cal I}{\cal K}}$ can take any positive value (for small $\delta>0$ take $a(t)=(1/\delta-1){\bf 1}_{[0,\delta]}(t)-{\bf 1}_{(\delta,1]}(t)$ or $a(t)={\bf 1}_{[1/2-\delta,1/2)}(t)-{\bf 1}_{[1/2,1/2+\delta]}(t)$). Also $e_{{\cal I}{\cal K}}$ can be arbitrarily close to $2e_{{\cal E}{\cal K}}$ 
(for small $\delta > 0$ take $A(t)=\sqrt{t(1-t)}\{t^{\delta}{\bf 1}_{[0,1/2]}(t)+(1-t)^{\delta}{\bf 1}_{(1/2,1]}(t)\}$).\\

\vspace{2mm}
\noindent 
{\bf Remark 8.} To give some insight into asymptotic levels of the tests considered in Theorems 1 - 3 and Remark 1, set $\theta_n=cn^{-q},$ where$ \;q\in(0,1/2),$ while $c$ is a positive constant. Recall that
the Kallenberg efficiency is characterized by levels $\alpha_n$ tending to 0 and asymptotic powers in $(0,1)$. According to (i) of Theorem 1 in Inglot et al. (2018), for any of the statistics, say ${\cal U}_n$, being compared to ${\cal K}_n$, it holds that $\log\alpha_n\sim -c_{\cal U}[b_{\cal U}(P_{\theta_n}^n)]^2$, where 
$c_{\cal U}[b_{\cal U}(P_{\theta_n}^n)]^2$ is the intermediate slope of ${\cal U}_n$.

For ${\cal G}_n$ and any $q\in (0,1/2)$ the allowable levels are of the form $\log\alpha_n \sim -\frac{c^2}{2}(\sup_{\kappa\leq t\leq 1-\kappa}|A^*(t)|)^2 \times n^{1-2q}$. 

For ${\cal E}_n$ take $\kappa_n\asymp n^{-\epsilon},\;\epsilon\in(0,1)$, and $A(t)$ satisfying (5.3). 
Then (7.2) holds for any $q\in(\epsilon/2,1/2)$ and the allowable levels take the form $\log\alpha_n \sim -\frac{c^2}{2}(\sup_t|A^*(t)|)^2 \times n^{1-2q}$.

For ${\cal M}_n$ take $A$ satisfying (5.6) with some $\varpi \in [0,1/2)$. Then (7.4) holds true for any $q\in(\varpi,1/2)$ and the allowable levels take the form $\log\alpha_n \sim-\log\left[c^2 n^{1-2q}(\sup_t|A^*(t)|)^2\right]\sim \log n^{2q-1}$.

For ${\cal I}_n$ the situation is much more regular. For any $q\in(0,1/2)$ in $\theta_n=cn^{-q}$ 
and $A$ satisfying (6.3) the statement 
(7.3) holds and the allowable levels take the form $\log\alpha_n\sim - c^2||A^*||_2^2 \times n^{1-2q}$. 

As for the asymptotic power under $F_n$ with the above $\theta_n$, we have the following situation, being a consequence of Lemma 1 of Appendix B in Inglot et al. (2018). In the case ${\cal I}_n$ any fixed asymptotic power from (0,1) is attainable by an appropriate choice of $w$ in (6.4). 
In contrast, 
for ${\cal G}_n,\; {\cal E}_n$ and ${\cal M}_n$ we do not show that asymptotic power exists and
we can only say that taking, in the present Theorems 1 and 2, any $w>0$ the resulting sequences of powers are bounded away from 0 and 1.\\

Though the above conclusions, contained in Remark 8, may look to be complicated and abstract, it turns out that, under standard circumstances, the value of the efficiency nicely helps to predict the empirical power of a test being compared to a benchmark. The reason for this is that on not very extreme tails of the test statistic, which are characteristic to the intermediate approach, the asymptotics work well for relatively small sample sizes. Hence, the approach gives good approximation for standard significance levels. A similar conclusion can also be found in Ermakov (2004), p. 624. Below, we
demonstrate to what extent, for selected statistics with non-zero intermediate efficiency with respect to ${\cal K}_n$, our results explain empirical powers under fixed levels and fixed alternatives.\\

\noindent
{\bf {\large 8. Simulation and efficiencies}}\\

\noindent
{\bf {8.1. Examples of departures from the standard Gaussian model }}\\

\noindent
We start with three simple classical situations related to detecting lack-of-fit to the standard normal distribution $N(0,1)$. To be specific, $F_0(x)=\Phi(x)$, and the alternatives are: $H_1(x;\mu)=\Phi(x-\mu),\;H_2(x;\sigma)=\Phi(x/\sigma)$ and $H_3(x;\mu,p)=(1-p)\Phi (x) + p \Phi (x - \mu).$ 
In all simulations here and in Section 9 we consider fixed alternatives. To clearly distinguish this case from the combination $(1-\vartheta_n)F_0+\vartheta_nF_1$, used in theoretical considerations, we use the notation $H_j,\;j=1,2,...$ for the fixed alternative. This is especially useful in Section 9, where $F_1$ itself corresponds to some mixtures.
For some simulated powers of ${\cal S}_n$ under the shift and scale models see Moscovich et al. (2016). The location-contaminated alternative $H_3(x;\mu,p)$ comes from the paper by Pearson et al. (1977). 
The alternative $H_3(x;\mu,p)$ was exploited for comparison of powers in Li and Siegmund (2015).
In recent years this model with $p=p_n,\;p_n \to 0$, and $\mu=\mu_n,\;\mu_n \to \infty$, has been popularized under the label ``sparse heterogeneous mixtures''; cf. Donoho and Jin (2004) and related papers.

After the transformation $\Phi(X_i),\;i=1,...,n$, these alternatives have some densities $h_j$ on $(0,1)$ which can always be written in the form $1+a^{[j]}(t)$, where $\int a^{[j]}(t) dt =0,\;j=1,2,3.$ Since we like to present $a^{[j]}$'s in our figures in some normalized form, we introduce the following parametrization.
By $||\cdot||_1$ we denote the $L_1$ norm on (0,1) with the Lebesgue measure, we put $\varphi= \Phi'$, $\theta^{[j]} = ||a^{[j]}||_1$ and $a_j = a^{[j]}/\theta^{[j]},\;j=1,2,3$. This yields the following alternative models:\\

\noindent
{$\mathbb{ M}_1 $:} $\displaystyle 
h_1(t;\mu)=1+\theta^{[1]} a_1(t;\mu),\;\;\;\mbox{with}\;\;\;a^{[1]}(t;\mu)=\frac{\varphi(\Phi^{-1}(t)-\mu)}{\varphi(\Phi^{-1}(t))} -1,\;\;\mu \in \mathbb{R},\; \mu \neq 0,$

\noindent
{$\mathbb{ M}_2 $:} $\displaystyle
h_2(t;\sigma)=1+\theta^{[2]} a_2(t;\sigma),\;\;\;\mbox{with}\;\;\;a^{[2]}(t;\sigma)=\frac{\varphi(\frac{1}{\sigma} \Phi^{-1}(t))}{\sigma \varphi(\Phi^{-1}(t))} -1,
\;\;\sigma \in \mathbb{R}_{+},\; \sigma \neq 1,$

\noindent
{$\mathbb{ M}_3 $:} $\displaystyle
h_3(t;p,\mu)=1+\theta^{[3]} a_3(t;p,\mu),\;\;\;\mbox{with}\;\;\;a^{[3]}(t;p,\mu)=p\Bigl\{\frac{\varphi(\Phi^{-1}(t)-\mu)}{\varphi(\Phi^{-1}(t))} -1\Bigr\},\\
\hspace{1cm} p \in (0,1),\; \mu \in \mathbb{R},\;\mu \neq 0.$
\\

The functions $a_1$ and $a_3$ are unbounded while $a_2$ is bounded for $\sigma \leq 1$ and unbounded otherwise. It holds that $a_j(u;\cdot) \in L_2(0,1),\;j=1,2,3$.
Set $A_j(t;\cdot)=\int_0^t a_j(u;\cdot) d u.$ We have
$A_1(t;\mu)=[\Phi(\Phi^{-1}(t) - \mu) -t]/\theta^{[1]},\;A_2(t;\sigma)=[\Phi(\frac{1}{\sigma}\Phi^{-1}(t))-t]/\theta^{[2]},\;A_3(t;p,\mu)=[pA_1(t;\mu)]/\theta^{[3]}=A_1(t;\mu).$ The last relation implies that the intermediate efficiency of the mixture does not depend on $p$. In contrast, the efficiency is influenced by a change of the ``direction'' of the noise in the mixture; i.e. $\Phi(x-\mu)$ in this particular case. More examples of mixtures are discussed in Section 9.

Similarly as in Section 5, given $A_j$, set 
$$
A_j^{*}(t)=\frac{A_j(t)}{\sqrt{t(1-t)}}.
$$
Note that for the functions $A_1$ and $A_3$ and all related parameters under consideration (5.6) holds with any $\varpi \in (0,1/2)$ and hence (5.3) and (6.3) hold, as well (cf. Remark 2).
For $A_2$, if $\sigma<1$ then (5.6) holds with $\varpi=0$; if $\sigma \in (1,\sqrt 2)$ then (5.6) holds with $\varpi \in [1-\sigma^{-2},1/2)$; if $\sigma=\sqrt{2}$ then (5.3) holds while (5.6) does not. For all $\sigma>0$ (6.3) is satisfied.\\

\noindent
{\bf 8.2. Alternatives from $\mathbb{ M}_1$, $\mathbb{ M}_2$ and $\mathbb{ M}_3$ satisfying (5.3), (5.6) and (6.3), corresponding efficiencies and simulated powers}\\

\noindent
We restrict our attention to ${\cal I}_n$, ${\cal M}_n$, ${\cal K}_n$, and two selected members of the class of statistics ${\cal E}_n={\cal E}_n(\kappa_n)$, indexed by $\kappa_n$ satisfying (ii) of Lemma 2. It is intuitively clear that using a relatively small parameter $\kappa_n$ can be profitable when under an alternative a considerable amount of a probability mass is shifted towards one or both tails, while a larger $\kappa_n$ is expected to be more useful in detecting centrally located changes. For an illustration we took
$$
{\cal E}_n^o={\cal E}_n^o(\kappa_n)\;\;\;\mbox{with}\;\;\;\kappa_n=\kappa_n^o=\frac{1}{2}n^{-1/2}
$$
and
$$
{\cal E}_n^{\star}={\cal E}_n^{\star}(\kappa_n)\;\;\;\mbox{with}\;\;\;\kappa_n=\kappa_n^{\star}=n^{-9/10}.
$$

In the simulation experiments the significance level was set to $\alpha=0.01$ and the number of MC runs for estimating sizes was $10^5$. Moreover, we used $10^4$ MC runs for estimating powers. The programs were written in C Sharp.

We have considered $\mathbb{M}_1$ with $\mu=0.15$, $\mathbb{M}_2$ with $\sigma=0.75$ and $\sigma=1.25$ and $\mathbb{M}_3$ with $p=0.05,\;\mu=2.00$. For all the cases the assumptions (5.3), (5.6) and (6.3) are satisfied. Hence our theoretical results on the intermediate efficiencies are applicable. 

The selected models, the corresponding efficiencies and the related empirical powers are presented in Figures 1 and 2. In the first row of the figures we display 
graphs of $a_j$ and $A^*_j$, $j=1,2,3$, and the corresponding values of $t_0$, $m_0$, where $t_0=\arg \max |A_j^*(t)|$ and $m_0=|A_j^*(t_0)|$. 

The middle rows show empirical powers of ${\cal E}_n^o$, ${\cal E}_n^{\star}$, ${\cal I}_n$, ${\cal M}_n$ and ${\cal K}_n$, against $n$. 

The bottom rows show the above power curves for sample sizes not exceeding the first value for which the empirical power of ${\cal E}_n^o$ attains the value in $[0.99,1]$. 
We additionally display here the values of the efficiencies $e_{{\cal E}{\cal K}}$ and $e_{{\cal I}{\cal K}}$. In all four cases $e_{{\cal E}{\cal K}} > 1$ as well as $e_{{\cal I}{\cal K}} > 1.$ In the last row
we also present the corresponding simulation results for ${\cal K}_{n \cdot e_{{\cal E}{\cal K}}}$ and ${\cal K}_{n \cdot e_{{\cal I}{\cal K}}}$ i.e. the empirical powers for the Kolmogorov-Smirnov test based on the corrected sample sizes ${n}\cdot{e_{{\cal E}{\cal K}}}$ and ${n}\cdot{e_{{\cal I}{\cal K}}}$, respectively. 
The zoom applied here allows to see well the way in which the corrected sample sizes influence the empirical powers of ${\cal K}_n$.

The results show that the finite sample interpretation of the intermediate efficiency indeed reflects very well the behavior of empirical powers of ${\cal K}_n$. 
For very large values of the efficiency $e_{{\cal E}{\cal K}}$ and relatively small sample sizes, as is the case for the model ${\mathbb M}_3$ in Figure 2, the empirical powers of ${\cal K}_{n \cdot e_{{\cal E}{\cal K}}}$ considerably overestimate the powers of ${\cal E}_n^o$ and ${\cal E}_n^{\star}$. However, it is hard to expect very accurate small sample results in such an extreme situation. In any case, the message is informative.
The results of simulations also indicate that the 0 efficiency of ${\cal M}_n$ with respect to ${\cal K}_n$ should not be surprising. Shapes of empirical powers of ${\cal M}_n$, as functions of $n$, are very different from those for ${\cal K}_n$. For the alternatives under consideration one needs a relatively huge number of observations to achieve a high power of the test based on ${\cal M}_n$. Similar pictures are expected to be valid for many other classical alternative distribution models.

\newpage

\begin{figure}[h!]
\centering

\includegraphics[scale=0.75]{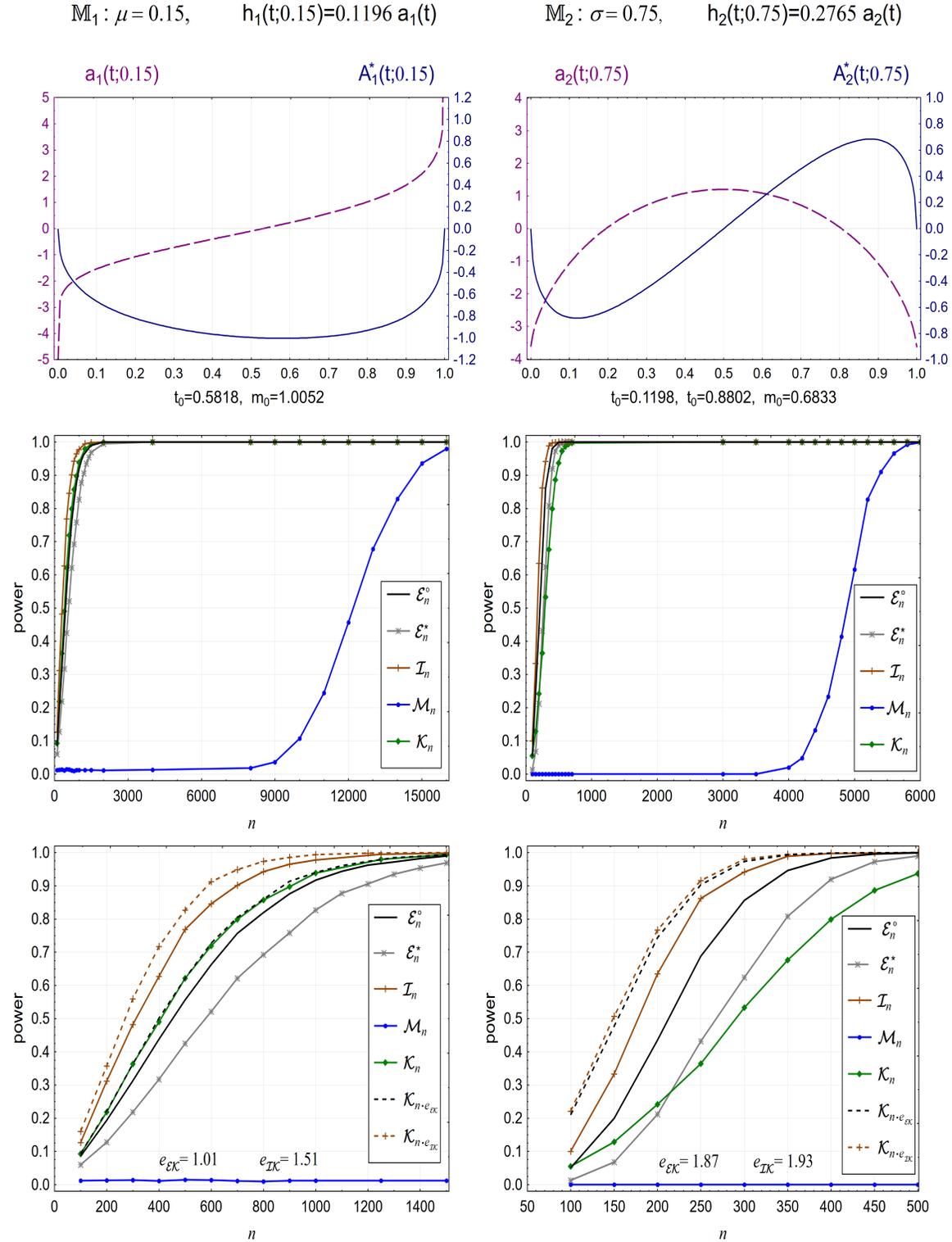}

\vspace{0.5cm}
\caption{{(\it Best viewed in color)} Alternatives from $\mathbb{M}_1,\;\mathbb{M}_2$ and empirical powers. {\it First row\;}: the functions $a_j$ - dashed line, and $A^*_j$ - solid line, 
and values of $t_0$, $m_0$, $j=1, 2$. {\it Second and third rows\;}: empirical powers (in the full range and zoomed) of ${\cal E}_n^o$, ${\cal E}_n^{\star}$, ${\cal I}_n$, ${\cal M}_n$, ${\cal K}_n$, ${\cal K}_{n \cdot e_{{\cal E}{\cal K}}}$, and ${\cal K}_{n \cdot e_{{\cal I}{\cal K}}}$. The third row also includes the corresponding efficiencies $e_{{\cal E}{\cal K}}$ and $e_{{\cal I}{\cal K}}$.}

\end{figure}
\newpage

\begin{figure}[h!]
\centering

\includegraphics[scale=0.75]{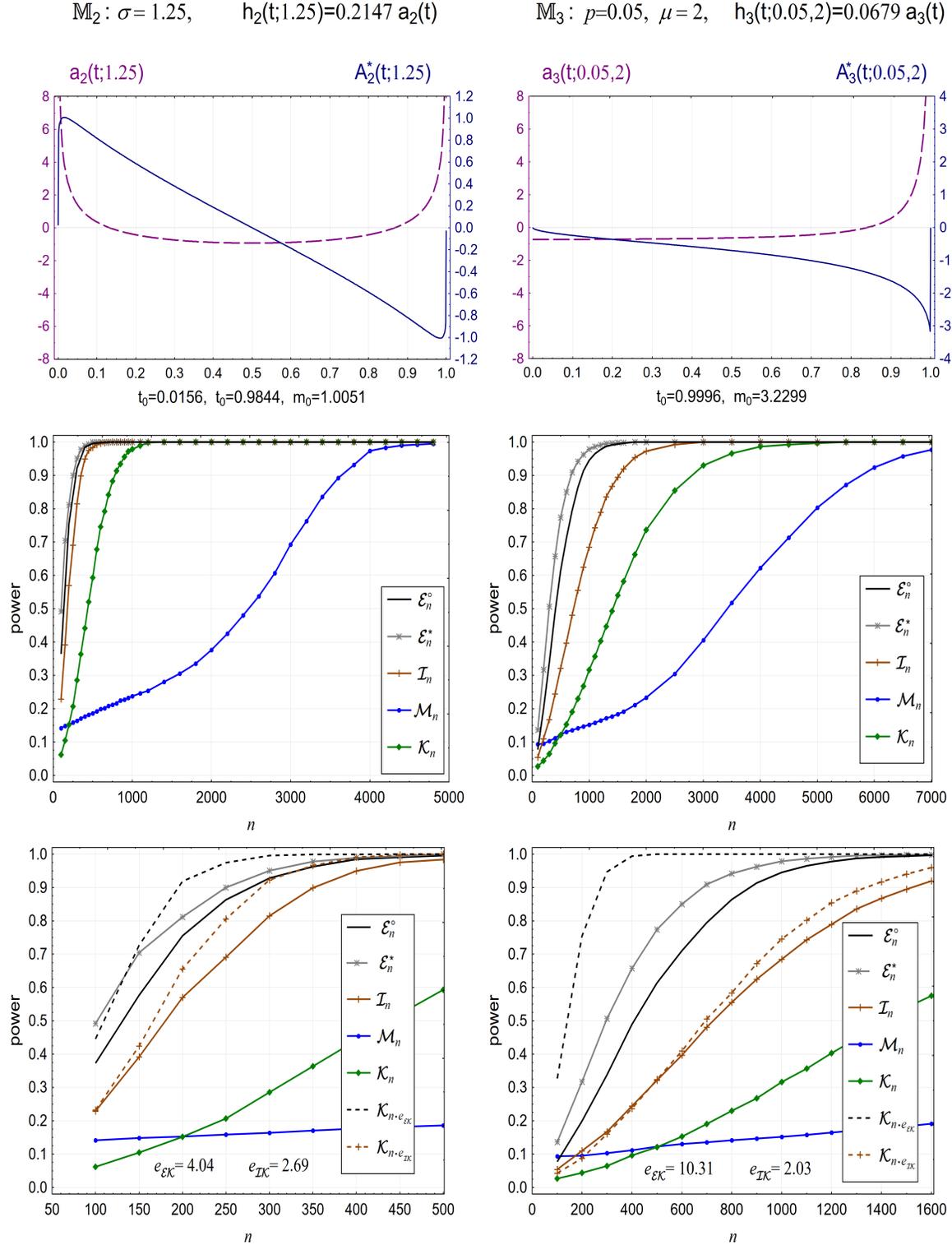}

\vspace{0.5cm}
\caption{{(\it Best viewed in color)} Alternatives from $\mathbb{M}_2$, $\mathbb{M}_3$ and empirical powers. {\it First row\;}: the functions $a_j$ - dashed line, and $A^*_j$ - solid line, 
 values of $t_0$, $m_0$, $j=2, 3$. {\it Second and third rows\;}: empirical powers (in the full range and zoomed) of ${\cal E}_n^o$, ${\cal E}_n^{\star}$, ${\cal I}_n$, ${\cal M}_n$, ${\cal K}_n$, ${\cal K}_{n \cdot e_{{\cal E}{\cal K}}}$, and ${\cal K}_{n \cdot e_{{\cal I}{\cal K}}}$. The third row also includes the corresponding efficiencies $e_{{\cal E}{\cal K}}$ and $e_{{\cal I}{\cal K}}$. }

\end{figure}

\noindent
{\bf {9. On the behavior of ${\cal E}_n$ and ${\cal K}_n$ when (5.3) is violated}}\\

\noindent
The above part of the paper gives some quite reliable insight into the behavior of powers of the Kolmogorov-Smirnov ${\cal K}_n$ test and the selected Eicker-Jaeschke statistics ${\cal E}_n$ and ${\cal M}_n$, in the case when the tails of an alternative are relatively light; i.e. the conditions (5.3) and (5.6) are satisfied. Under these conditions $e_{{\cal E}{\cal K}} \geq 1$ and $e_{{\cal M}{\cal K}}=0$, respectively. From previous developments it follows that one should expect much worse power behavior of ${\cal K}_n$ in the case of alternatives obeying relatively heavy tails. We shall study this question in the present section by contrasting the behavior of ${\cal K}_n$ with ${\cal E}_n$, in the case when the condition (5.3) is violated. Since we are aware of an extension of Theorem 1 in this case, we are able to calculate only a so-called weak variant of the intermediate efficiency. Let us denote it by $\hat e_{{\cal E}{\cal K}}$.
This weak variant is defined as a limit of the ratio of the slopes, as $n$ tends to infinity. The difference between $\hat e_{{\cal E}{\cal K}}$ and $\ e_{{\cal E}{\cal K}}$ resembles to some extent the difference between the approximate and the exact Bahadur efficiency. 
The weak variant of the intermediate efficiency was already studied in Ivchenko and Mirakhmedov (1995), and Inglot (1999).

To calculate $\hat e_{{\cal E}{\cal K}}$ for a local sequence of alternatives $F_n(t)=t+\theta_n A(t)$, when (5.3) is violated, set 
$$
m_n= \sup_{\kappa_n \leq t \leq 1-\kappa_n} \frac{|A(t)|}{\sqrt{t(1-t)}}
\eqno(9.1)
$$
and denote by $t_n$ any point at which the supremum in (9.1) is attained. \\

\noindent
{\bf Lemma 7.} {\it Suppose that $m_n \to \infty$ and $t_n \to 0$ or $t_n \to 1$, as $ n \to \infty$. Assume that $\lim \inf_{n \to \infty} n\kappa_n/\log ^2 n > 0$, 
$\lim_{n \to \infty} \log \kappa_n /\log n < 0, \;n\theta_n^2/\log \log (1/\kappa_n) \to \infty$, and $\theta_n^2 m_n^2/\kappa_n \to 0$. Then one gets }
$$
\lim_{n \to \infty} P_{\theta_n}^n \Bigl(\Big|\frac{{\cal E}_n}{\sqrt n \theta_n m_n} -1\Big| \leq \epsilon\Bigr)=1\;\;\;\mbox{for every}\;\;\;\epsilon > 0.
\eqno(9.2)
$$
{\it Hence, the intermediate slope ${c_{\cal E}[b_{\cal E}(P_{\theta_n}^n)]^2 }$ of ${\cal E}_n$ under $\{P_{\theta_n}\}$ has the form $n\theta_n^2 m_n^2 /2$}.\\

\noindent
{\bf Corollary 1}. {\it Under the assumptions of Lemma 7 it holds that}
$$
\hat e_{{\cal E}{\cal K}} = \lim_{n \to \infty} \frac{c_{\cal E}[b_{\cal E}(P_{\theta_n}^n)]^2 }{c_{\cal K}[b_{\cal K}(P_{\theta_n}^n)]^2 }=
\lim_{n \to \infty}\frac{m_n^2}{4||A||_{\infty}^2} = +\infty.
\eqno(9.3)
$$
\\

The relation (9.3) suggests that perhaps the intermediate efficiency $e_{{\cal E}{\cal K}}$ of ${\cal E}_n$ with respect to ${\cal K}_n$ equals $+ \infty$, as well. However, verifying this would require non-trivial investigations of the question on non-degeneracy of the asymptotic power of ${\cal E}_n$ under the above described local alternatives. This is a challenging open question. Note that non-degenerate asymptotic power of ${\cal E}_n$ is needed to have the interpretation of the intermediate efficiency in terms of the limiting ratio of appropriate sample sizes; cf. Theorem 1 in Inglot et al. (2018).

We show below that even this weak variant $\hat e_{{\cal E}{\cal K}}$ of the efficiency gives a right indication on an empirical power behavior of ${\cal E}_n$ and ${\cal K}_n$, when (5.3) fails. 

We shall study an empirical behavior of ${\cal E}_n^o$, ${\cal E}_n^{\star}$, ${\cal K}_n$, as well as ${\cal M}_n$ and ${\cal I}_n$ under the following alternative models

\vspace{0.2cm}
\noindent
${\mathbb M}_4 : H_4(t;\beta, \pi) = \{\pi^{(\beta-1)/\beta} {t^{1/\beta}}\}{\bf 1}_{[0,\pi)}(t) + t {\bf 1}_{[\pi,1-\pi]}(t) + \{1-\pi^{(\beta -1)/\beta}(1-t)^{1/\beta}\}{\bf 1}_{(1-\pi,1]}(t)$, where $\beta >0,\;\pi \in [0,0.5]$, and $t \in [0,1]$,

\vspace{0.2cm}
\noindent 
${\mathbb M}_5 : H_5(x;\delta,p)=(1-p)\Phi(x) +p \Lambda(x;\delta),\; \delta >0, \;p \in [0,1],\; x \in \mathbb{R},$ where $\Lambda(x;\delta)=[\Phi(x)]^{\delta}$ is the Lehmann distribution; cf. Remark 3,

\vspace{0.2cm}
\noindent
${\mathbb M}_6 : H_6(x;\gamma,p)=(1-p)\Phi(x) +p \Sigma(x;\gamma),\; \gamma > 0, \;p \in [0,1],\; x \in \mathbb{R},$ where $\Sigma(x;\gamma)$ is the symmetric Subbotin distribution function obeying the density $C_{\gamma} \exp\{-|x|^{\gamma}/\gamma\},\; x \in \mathbb{R},$ 

\vspace{0.2cm}
\noindent
${\mathbb M}_7 : H_7(x;\zeta,p)=(1-p)\Phi(x) +p \Pi(x;\zeta),\; \zeta > 0, \;p \in [0,1],\; x \in \mathbb{R},$ where $\Pi(x;\zeta)$ is the distribution function of the symmetric Pareto distribution with the parameter $\zeta$; cf. Remark 4.

\vspace{0.2cm}
The model ${\mathbb M}_4$ comes from Mason and Schuenemeyer (1983). If $\beta \in (0,1)$ then $H_4(t;\beta,\pi)$ has lighter tails than the uniform (0,1) distribution, say $U(0,1)$. When $\beta > 1$ then $H_4(t;\beta,\pi)$ has heavier lower and upper tails than $U(0,1)$. For ${\mathbb M}_4$ the condition (5.3) does not hold if $\beta \geq 2$. ${\mathbb M}_4$ defines alternatives with an allocation of the probability mass only on the tails.

${\mathbb M}_5$ - ${\mathbb M}_7$ were chosen as mixtures. Detection of mixtures is of vital interest. Lehmann's model, used in ${\mathbb M}_5$, is popular in the statistical literature. The Subbotin distribution is discussed in Donoho and Jin (2004). The mixture ${\mathbb M}_7$ has been inspired by Jin et al. (2005), where an additive model with disturbances with algebraically decreasing tails was considered. For ${\mathbb M}_5$ with $\delta >0$, ${\mathbb M}_6$ with $\gamma \in (0,2)$, and ${\mathbb M}_7$ with $\zeta > 0$ the condition (5.3) does not hold. 

Each of the models ${\mathbb M}_j,\; j=4,...,7,$ can be equivalently rewritten in the form $1+\theta^{[j]} a_j(t;\cdot). $ 
The functions $a_j,\; j=4,6,7$, are symmetrical with respect to 1/2 and unbounded at 0 and 1 while $a_5$ is unbounded at 0. For $t$ close to 0 the functions $a_4,..., a_7$ behave like: $ t^{(1-\beta)/\beta},\; t^{\delta-1},\; t^{-1}\exp\{-\frac{1}{\gamma}[\log(1/t^2\log(1/t^2))]^{\gamma/2}-\frac{1}{2}[\log\log(1/t)]\}$, $\;t^{-1}[\log(1/t)]^{-1-\zeta/2},$ respectively. Note also that $a_4,...,a_7$ do not belong to $L_2(0,1)$ for $\beta \geq 2, \delta \leq1/2, \gamma <2, \zeta >0 $, accordingly.

In Figure 3 we plot empirical powers of the considered tests, under $\alpha=0.01$ and some selected $n$ and $p$, against the parameters $\pi,\;\beta,\;\delta,\;\gamma,$ and $\zeta$ of the considered models. The outcomes show that, when (5.3) is violated, empirical behavior of ${\cal K}_n$ is very poor and resembles the behavior of ${\cal M}_n$
in previous figures. In contrast, now ${\cal M}_n$ does very well. Obviously, the imposed lack of (5.3) implies the violation of (5.6), as well.
Moreover, except for the cases when a very large amount of probability mass is shifted to the ends of $(0,1)$, ${\cal E}_n^{\star}$ also works very well. In all situations shown in Figure 3 the variant ${\cal E}_n^{\star}$ dominates ${\cal E}_n^o$ considerably. The empirical behavior of ${\cal I}_n$ is not impressive in comparison to ${\cal M}_n$ and ${\cal E}_n^{\star}$.

It should be emphasized that we have not conducted an extensive search for $\kappa_n^o$ and $\kappa_n^{\star}$ defining ${\cal E}_n^{o}$ and ${\cal E}_n^{\star}$. We simply took the two candidates which satisfy the assumption (ii) of Lemma 2, i.e. $\kappa_n$ satisfying $\liminf_{n} n \kappa_n /\log^2 n > 0$. In spite of this, from the outcomes in Figures 1 - 3, it can be seen that ${\cal E}_n^{\star}$ is a reasonably well balanced solution. At any rate, some search for a data-driven choice of the smoothing parameter $\kappa_n$ would be very welcome.
\\

\begin{figure}[h!]
\centering

\includegraphics[scale=0.75]{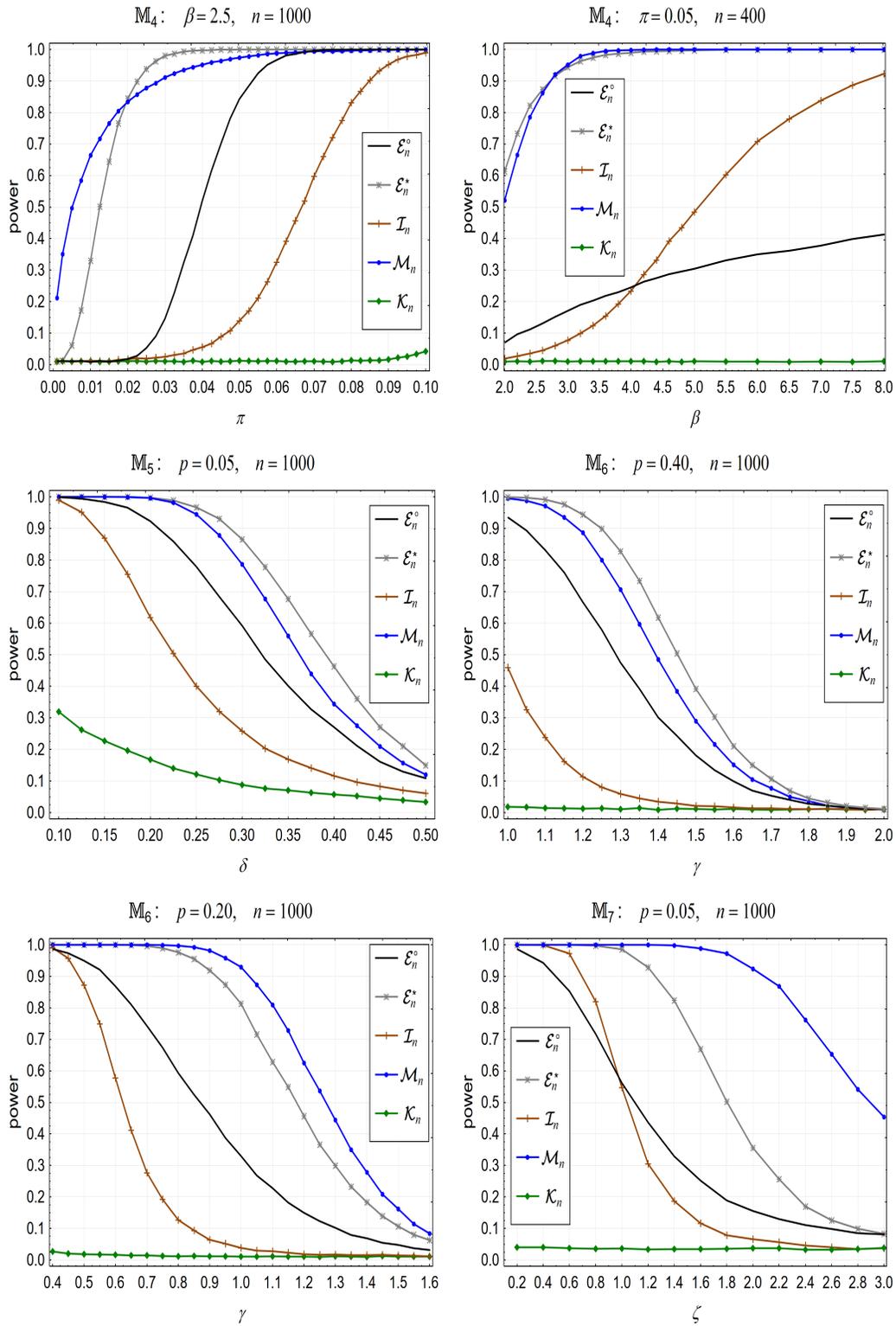}

\vspace{0.5cm}
\caption{{(\it Best viewed in color)} Empirical powers of ${\cal E}_n^o$, ${\cal E}_n^{\star}$, ${\cal K}_n$, ${\cal I}_n$, and ${\cal M}_n$ under alternatives $\mathbb{M}_4 - \mathbb{M}_7$ for selected parameters and sample sizes. }

\end{figure}
\newpage

\noindent
{{\bf \large 10. Discussion} \\

\noindent 
The present paper illustrates the advantages of using the pathwise variant of the Kallenberg efficiency to study goodness-of-fit to a completely known continuous distribution function. In Inglot et al. (2018) the paths were defined as mixtures of a big fraction of the null distribution and a small fraction of an alternative one. Consequently, we consider $(1-\vartheta_n)F_0(x) + \vartheta_n F_1(x)$, where $F_0$ is the null distribution, $F_1$ represents the alternative, and $\vartheta_n \to 0$ as $n \to \infty$. For convenience, in this paper we have transformed the observations to (0,1) via $F_0$, cf. (2.1), but it is not essential to the interpretation of the results. Moreover, to increase the readability of the results, we introduced (2.2). Anyway, in essence the pathwise variant of the efficiency evaluates the quality of tests by measuring their ability to detect (local) mixtures.
On the other hand, the mixtures define ``directions'' along which we approach the null model and, as a rule, the corresponding results on the efficiency are valid for many ``directions''. Moreover, in the intermediate approach $\vartheta_n$ decreases relatively slowly. The above implies that the resulting, asymptotic in nature, expression for the efficiency gives reliable results on empirical powers under fixed alternatives which are not necessarily mixtures, fixed sample sizes, and standard significance levels. 

At first glance, our approach resembles detecting mixtures under the dense regime; cf. Cai et al. (2011) for the terminology and an insightful introduction to the problem. However, we are focused on a goodness-of-fit context and our goal is not to study if and when a procedure can detect or fail to detect a given mixture, but we would like to investigate how well a selected test can distinguish some classes of alternatives from the null model. Therefore, in contrast to the signal detection approach, we insist on having the error of the second kind in $(0,1)$. Moreover, the distribution function $F_1$ is fixed, independent on $n$. So, our setting differs from the typical approach in studies of detectable and undetectable regions, originated by Ingster (1997) and extensively developed in recent years; cf. Ditzhaus (2018) for the most general setting and historical details. Also, the outcomes of both approaches are qualitatively different. A typical feature of Ingster's approach is that whole big classes of tests have the same detection boundaries; cf. Jager and Wellner (2007), and Ditzhaus (2018) for an illustration. In contrast, the Kallenberg efficiency allows for catching some subtle differences between test statistics. It seems that some further investigations on this approach could result in better understanding advantages and limitations of popular classes of modern goodness-of-fit statistics. In particular, some more work on the asymptotic distribution of test statistics under the regime $\vartheta_n \to 0$ and $n \vartheta_n^2 \to \infty$ is necessary. Moreover, moderate deviations for the whole classes of test statistics, which were recently considered, should be developed. As illustrated by our analysis of ${\cal I}_n$ and related discussion, for sufficiently smooth functionals of the weighted empirical process deriving the intermediate efficiency is relatively easy. Sup-type functionals are less regular and more difficult to handle. Anyway, in our opinion, the present paper shows that such work is worthy of further consideration. In particular, it would be interesting to close our investigations on ${\cal M}_n$ and ${\cal E}_n$ by showing if and when their intermediate efficiencies with respect to ${\cal K}_n$ exist in the situation when (5.6) and (5.3), respectively, are violated.\\


\noindent
{\large {\bf Appendix: Proofs}}\\

\noindent
{\bf A.1. Proof of Lemma 2}\\

\noindent
Let $U_1,...,U_n$ be independent uniform (0,1) random variables and let $U_{(1)} \leq ... \leq U_{(n)}$ denote their order statistics.

(i) Let $i_n=\lfloor 3n w_n\sqrt{\kappa_n}\rfloor$. Then, by the assumption $w_n/\sqrt{\kappa_n}\to \infty$, we have for sufficiently large $n$
$$ P_0^n({\cal E}_n\geq \sqrt{n}w_n)\geq P_0^n\left(\max_{\{i:U_{(i)}\in [\kappa_n,1-\kappa_n]\}}\frac{|U_{(i)}-i/n|}{\sqrt{U_{(i)}(1-U_{(i)})}}\geq w_n\right)\hspace*{3cm}$$
$$\hspace*{0.4cm}\geq P_0^n\left(|U_{(i_n)}-i_n/n|\geq w_n \sqrt{U_{(i_n)}},\;\, \kappa_n\leq U_{(i_n)}\leq 1-\kappa_n\right)$$
$$\hspace*{1.5cm}\geq P_0^n\left(U_{(i_n)}-i_n/n\leq -w_n \sqrt{U_{(i_n)}}\right)-P_0^n(U_{(i_n)}< \kappa_n)-P_0^n(U_{(i_n)}>1-\kappa_n)$$
$$\geq P_0^n\left(U_{(i_n)}\leq i_n^2/4n^2w_n^2\right)-P_0^n(U_{(i_n)}<\kappa_n)-P_0^n(U_{(i_n)}>1-\kappa_n).\eqno (A.1)$$
Since $j!\geq j^{j}e^{-j}$ for all $j\geq 1$ then 
$$P_0^n(U_{(i)}\leq u)=\sum_{j=i}^n\left(\begin{array}{c} n\\j\end{array}\right)u^j(1-u)^{n-j}\leq\sum_{j=i}^n\left(\frac{enu}{j}\right)^j\leq\sum_{j=i}^n\left(\frac{enu}{i}\right)^j.$$
Hence and from the relation $en\kappa_n/i_n<1/2$ for sufficiently large $n$ we get
$$ P_0^n(U_{(i_n)}\leq \kappa_n)\leq\left(\frac{en\kappa_n}{i_n}\right)^{i_n}\frac{i_n}{i_n-en\kappa_n}\leq 2\left(\frac{e}{3}\frac{\sqrt{\kappa_n}}{w_n}\right)^{i_n}. \eqno (A.2)$$
Moreover, since $i_n/n\to 0$, then for sufficiently large $n$ it holds $P_0^n(U_{(i_n)}>1-\kappa_n)\leq P_0^n(U_{(i_n)}<\kappa_n)$. On the other hand by $j!\leq j^{j+1}e^{-j}$ being true for $j\geq 7$ we have
$$ P_0^n(U_{(j)}\leq u)\geq \left(\begin{array}{c} n\\j\end{array}\right)u^j(1-u)^{n-j}\geq \frac{[(n-j)eu]^j}{j^{j+1}}(1-u)^{n}.$$
As $1-i_n/n> 2/3$ for sufficiently large $n$, the above inequality and the definition of $i_n$ imply for sufficiently large $n$ 
$$P_0^n\left(U_{(i_n)}\leq i_n^2/4n^2w_n^2\right)
\geq\frac{1}{i_n}\left(1-\frac{i_n^2}{4n^2w_n^2}\right)^n\left(\frac{(n-i_n)ei_n}{4n^2w_n^2}\right)^{i_n}\hspace*{3cm}$$
$$\hspace*{3cm}\geq\frac{1}{i_n}(1-9\kappa_n/4)^n\left(\frac{e}{2}\frac{\sqrt{\kappa_n}}{w_n}\right)^{i_n}\geq \frac{1}{i_n}e^{-3n\kappa_n}\left(\frac{e}{2}\frac{\sqrt{\kappa_n}}{w_n}\right)^{i_n}.\eqno (A.3)$$
Combining (A.1), (A.2) and (A.3), again by the definition of $i_n$ and the assumption $w_n/\sqrt{\kappa_n}\to\infty$, we obtain for sufficiently large $n$
$$P_0^n({\cal E}_n\geq \sqrt{n}w_n)\geq \frac{1}{i_n}e^{-3n\kappa_n}\left(\frac{e}{2}\frac{\sqrt{\kappa_n}}{w_n}\right)^{i_n}-4\left(\frac{e}{3}\frac{\sqrt{\kappa_n}}{w_n}\right)^{i_n}=\frac{1}{i_n}e^{-3n\kappa_n}\left(\frac{e}{2}\frac{\sqrt{\kappa_n}}{w_n}\right)^{i_n}(1+o(1)). \eqno (A.4)$$
Imposing the logarithm in (A.4), dividing by $-nw_n^2$, and using again the assumption $w_n/\sqrt{\kappa_n}\to\infty$ we get
$$-\frac{1}{nw_n^2} \log P_0^n({\cal E}_n\geq \sqrt{n}w_n)\leq-\frac{3\sqrt{\kappa_n}}{w_n}\log \frac{e}{2}\frac{\sqrt{\kappa_n}}{w_n}+3\frac{\kappa_n}{w_n^2}+\frac{\log i_n}{nw_n^2}+o(1)\to 0$$
and the proof is complete.\hfill $\Box$

(ii) Let $u_n(t)$ be the uniform empirical process and denote
$$ Z_n=\sup_{[\kappa_n,1-\kappa_n]} \frac{|B(t)|}{\sqrt{t(1-t)}}.$$

Since $nw_n^2\to\infty$, then the assumption on $\kappa_n$ implies $n^2w_n^2\kappa_n/\log^2n\to\infty$. Let $\varepsilon_n>0,\;\varepsilon_n\to 0$, be such that $w_n^2/(\kappa_n\varepsilon_n^2)\to 0$ and $n^2w_n^2\kappa_n\varepsilon_n^2/\log^2n\to\infty$. Then for any fixed $c\in(0,1)$ and sufficiently large $n$ we have
$$P_0^n({\cal E}_n\geq \sqrt n w_n)\geq Pr(Z_n\geq (1+\varepsilon_n)\sqrt n w_n)-Pr(\sup_{(0,1)}|u_n(t)-B_n(t)|\geq c\varepsilon_n\sqrt{n \kappa_n}w_n)\eqno (A.5)$$
and
$$P_0^n({\cal E}_n\geq\sqrt n w_n)\leq Pr(Z_n\geq (1-\varepsilon_n)\sqrt n w_n)+Pr(\sup_{(0,1)}|u_n(t)-B_n(t)|\geq c\varepsilon_n\sqrt{n\kappa_n}w_n).\eqno(A.6)$$
Moreover, from KMT inequality we have
$$Pr(\sup_{(0,1)}|u_n(t)-B_n(t)|\geq c\varepsilon_n\sqrt{n\kappa_n}w_n)\leq L\exp\{-lc\varepsilon_n\sqrt{\kappa_n}nw_n+lC\log n\}, \eqno(A.7)$$
where $l,L,C$ are universal positive constants. 

If we shall show that for any $w_n\to 0$ and such that $nw_n^2\to\infty$ it holds
$$-\frac{1}{nw_n^2}\log Pr\left(Z_n\geq\sqrt n w_n\right)\to \frac{1}{2} \eqno(A.8)$$
then by the choice of $\varepsilon_n$, the first component in the exponent on the right hand side of (A.7) dominates the second one and simultaneously the first component on the right side of (A.5) and (A.6) dominates the second one and (4.4) follows from (A.8).

To prove (A.8) recall that from the Darling-Erd\"os theorem (cf. Cs\"org\H{o} and Horvath, 1993, pp. 257-258) it follows that 
$$ Pr(\sqrt{a_n}Z_n-a_n-\frac{1}{2}\log(a_n/(2\pi))\leq y)\to \exp\{-2e^{-y}\},$$
where $a_n=2\log\log(1/\kappa_n-1)$. Denote by $\mu_n$ the median of $Z_n$. Then from the above relation $\sqrt{a_n}\mu_n-a_n-(\log(a_n/(2\pi))/2\to \mu$, where $\mu=\log(2/\log 2)$ is the median of the limiting distribution. Hence $\mu_n=\sqrt{a_n}+o(1)$, and $\mu_n$ tends to infinity. By a straightforward application of the Borell inequality for $Z_n$ (see e.g. van der Vaart and Wellner, 2000, p. 438) we get for every $n$ and $y>0$
$$ Pr(|Z_n-\mu_n|\geq y)\leq \exp\{-y^2/2\}.$$
Since by the assumption it follows $nw_n^2/\mu_n^2\to\infty$ then inserting $y=\sqrt n w_n-\mu_n$ into the last inequality we get 
$$-\limsup_{n\to\infty}\frac{1}{nw_n^2}\log Pr(Z_n\geq\sqrt n w_n)\geq \frac{1}{2}.$$
On the other hand, for any $\epsilon\in (0,1/2)$ and sufficiently large $n$ we have $\kappa_n<\epsilon$ and consequently
$$-\liminf_{n\to\infty}\frac{1}{nw_n^2}\log Pr(Z_n\geq\sqrt n w_n)\leq -\liminf_{n\to\infty}\frac{1}{nw_n^2}\log Pr(\sup_{t\in[\epsilon,1-\epsilon]}\frac{|B(t)|}{\sqrt{t(1-t)}}\geq\sqrt n w_n)= \frac{1}{2}.$$
The last two relations complete the proof of (A.8).\hfill $\Box$\\

\noindent
{\bf A.2. Proof of Lemma 4 (i)}\\

\noindent
The proof goes along the lines of that of Lemma 1 in Mason (1985). For any $n\geq 1$ the function $h(y)=y+w_ny^{\tau}-1/n$ is increasing on $(0,\infty)$ and
$h((2nw_n)^{-1/\tau})<0$ due to $nw_n^{1/(1-\tau)}\geq nw_n^2\to\infty$. This gives the following estimate
$$ P_0^n({\cal C}_n\geq \sqrt{n}w_n)\geq P_0^n\left(U_{(1)}-\frac{1}{n}\leq -w_n U_{(1)}^{\tau}\right)$$
$$\geq P_0^n\left(U_{(1)}\leq \frac{1}{(2nw_n)^{1/\tau}}\right)=1-\left(1-\frac{1}{(2nw_n)^{1/\tau}}\right)^n.$$
By the inequality $1-(1-y)^n>ny/e$ holding for $y<1/n$ we infer that for some positive $c$
$$P_0^n({\cal C}_n\geq \sqrt{n}w_n)\geq c n^{1-1/\tau}w_n^{-1/\tau}. \eqno(A.9)$$
Taking logarithms of both sides of (A.9), dividing by $-nw_n^2$ and using the assumption $nw_n^2/\log n\to\infty$ we get
$$-\frac{1}{nw_n^2}\log P_0^n({\cal C}_n\geq \sqrt{n}w_n)\leq (\frac{1}{2\tau}-1)\frac{\log n}{nw_n^2}+\frac{1}{2\tau}\frac{\log nw_n^2}{nw_n^2}-\frac{\log c}{nw_n^2}\to 0 \eqno (A.10)$$
which completes the proof. \hfill $\Box$\\

\noindent
{\bf A.3. Proof of Lemma 3 (i)}\\

\noindent 
Observe that for $\tau=1/2$ the above proof is valid. The only difference is that in (A.10) the first component on the right hand side vanishes and the assumption $nw_n^2/\log n\to\infty$ becomes superfluous. So, Lemma 3 (i) holds true. \hfill $\Box$\\

\noindent
{\bf A.4. Proof of Theorem 1}\\

\noindent 
Let $u_n(t),\;t \in (0,1)$, be the uniform empirical process and set $v(t)=\sqrt{t(1-t)}$. By (5.3) there exists $t_0\in (0,1)$ such that $\displaystyle |A^*(t_0)|=\sup_{(0,1)}|A^*(t)|=m_0$.

{\bf (i)} It holds that
$$P_{\theta_n}^n({\cal E}_n-b_{\cal E}(P_{\theta_n}^n)\leq w)
=P_{\theta_n}^n\left({\cal E}_n\leq w+\sqrt n\theta_n\frac{|A(t_0)|}{\sqrt{t_0(1-t_0)}}\right)
$$
$$
\leq Pr\left(\frac{|u_n(F_n(t_0))+\sqrt n\theta_nA(t_0)|}{\sqrt{t_0(1-t_0)}}\leq w+\sqrt n\theta_n\frac{|A(t_0)|}{\sqrt{t_0(1-t_0))}}\right).
\eqno(A.11)
$$
When $A(t_0)>0$ then (A.11) is majorized by ${Pr}\left(u_n(F_n(t_0))/\sqrt{t_0(1-t_0)}\leq w\right)$ converging to ${\Phi}(w)$. 
If $A(t_0)<0$ then the majorant $ Pr\left(u_n(F_n(t_0))/\sqrt{t_0(1-t_0)}\geq -w\right)$ of (A.11) converges to ${\Phi}(w)$ as well. This proves (i).\\

{\bf (ii)}. The key step is to show that for $\delta=\delta_{\cal E}(A)$ appearing in (5.4)
$$ Pr \left(\sup_{[\kappa_n,1-\kappa_n]}\frac{|u_n(F_n(t))+\sqrt n\theta_nA(t)|}{v(t)}>\sup_{[\delta,1-\delta]}\frac{|u_n(F_n(t))+\sqrt n\theta_nA(t)|}{v(t)}\right)=o(1).
\eqno(A.12)
$$
Indeed, having (A.12), for positive $w$ the triangle inequality and (5.4) imply that
$$P_{\theta_n}^n({\cal E}_n-b_{\cal E}(P_{\theta_n}^n)\leq w)={Pr}\left(\sup_{[\kappa_n,1-\kappa_n]}\frac{|u_n(F_n(t))+\sqrt n\theta_nA(t)|}{v(t)}-\sup_{(0,1)}\frac{|\sqrt n\theta_nA(t)|}{v(t)}\leq w\right)$$
$$\geq { Pr}\left(\sup_{[\delta,1-\delta]}\frac{|u_n(F_n(t))+\sqrt n\theta_nA(t)|}{v(t)}-
\sup_{[\delta,1-\delta]}\frac{|\sqrt n\theta_nA(t)|}{v(t)}\leq w\right)+o(1)$$
$$\geq { Pr}\left(\sup_{[\delta,1-\delta]}\frac{|u_n(F_n(t))|}{v(t)}\leq w\right)+o(1).$$

Since $u_n \circ F_n$ converges in distribution to a Brownian bridge, then (ii) follows. 

Now, by the definitions of $m_0$ and $\delta$, using the triangle inequality we infer that
$$ { Pr}\left(\sup_{[\kappa_n,1-\kappa_n]}\frac{|u_n(F_n(t))+\sqrt n\theta_nA(t)|}{v(t)}>\sup_{[\delta,1-\delta]}\frac{|u_n(F_n(t))+\sqrt n\theta_nA(t)|}{v(t)}\right)$$
$$={ Pr}\left(\sup_{[\kappa_n,1-\kappa_n]\setminus[\delta,1-\delta]}\frac{|u_n(F_n(t))+\sqrt n\theta_nA(t)|}{v(t)}>\sup_{[\delta,1-\delta]}\frac{|u_n(F_n(t))+\sqrt n\theta_nA(t)|}{v(t)}\right)$$
$$\leq {Pr}\left(\sup_{[\kappa_n,1-\kappa_n]\setminus[\delta,1-\delta]}\frac{|u_n(F_n(t))|}{v(t)}+\frac{m_0}{2}\sqrt n\theta_n>m_0\sqrt n\theta_n-\sup_{[\delta,1-\delta]}\frac{|u_n(F_n(t))|}{v(t)}\right)
$$
$$
\leq { Pr}\left(\sup_{[\kappa_n,1-\kappa_n]}\frac{|u_n(F_n(t))|}{v(t)}>\frac{m_0}{4}\sqrt n\theta_n\right).\eqno (A.13)$$

For $t\in (0,1)$ we have
$$0 \leq \frac{F_n(t)(1-F_n(t))}{t(1-t)}= 1 + \theta_n\frac{1-2t}{t(1-t)}A(t)-\theta_n^2\frac{A^2(t)}{t(1-t)}\leq 1+\theta_n\frac{|1-2t|}{\sqrt{t(1-t)}}|A^*(t)|. 
\eqno (A.14)$$
So, by (5.3) and the assumption $\theta_n^2/\kappa_n\to 0$ for $t\in[\kappa_n,1-\kappa_n]$ and sufficiently large $n$ the right hand side of (A.14) can be estimated by
$$
1+\frac{\theta_n}{\sqrt{\kappa_n}}m_0\leq 2.$$

Hence, for $t\in[\kappa_n,1-\kappa_n]$ and $n$ sufficiently large we have
$$
\frac{|u_n(F_n(t))|}{v(t)} \leq 2\frac{|u_n(F_n(t))|}{\sqrt{F_n(t)(1-F_n(t))}},
$$
and the right hand side in (A.13) is majorized by
$$Pr\Bigl(\sup_{(0,1)} \frac{|u_n(t)|}{v(t)} > \frac{m_0}{8}\sqrt n \theta_n \Bigr)$$
which, in view of the assumption $(n \theta_n^2)/\log \log n \to \infty$, as $n \to \infty$, and an application of the main result of Mason (1985), tends to 0. This concludes the proof of (A.12). \hfill $\Box$\\

\noindent
{\bf A.5. Proof of Theorem 2}\\

\noindent 
As previously, let $u_n(t),\;t \in (0,1)$, be the uniform empirical process and set $v(t)=\sqrt{t(1-t)}$. By (5.6) there exists $t_0\in (0,1)$ such that $\displaystyle |A^*(t_0)|=\sup_{(0,1)} |A^*(t)|=m_0$. We can write
$$ 
{\cal S}_n=\sqrt{n}\sup_{(0,1)}\frac{|\hat{F}_n(t)-t|}{v(t)}\stackrel{ D}{=}\sup_{(0,1)}\frac{|u_n(F_n(t))+\sqrt n\theta_nA(t)|}{v(t)}.
$$

\noindent
{\bf Proof of (i)}. Since
$$P_{\theta_n}^n({\cal S}_n-\sqrt n\theta_nm_0\leqslant w)
=P_{\theta_n}^n\left({\cal S}_n\leqslant w+\sqrt n\theta_n\frac{|A(t_0)|}{v(t_0)}\right)
$$
$$
\leqslant Pr\left(\frac{|u_n(F_n(t_0))+\sqrt n\theta_nA(t_0)|}{v(t_0)}\leqslant w+\sqrt n\theta_n\frac{|A(t_0)|}{v(t_0)}\right),
$$
then we proceed exactly in the same way as in the proof of (i) in Theorem 1.\\

\noindent
{\bf Proof of (ii)}. The key step is to show that for $\delta=\delta_{\cal M}(A)$ appearing in (5.7)
$$ Pr \left(\sup_{(0,1)}\frac{|u_n(F_n(t))+\sqrt n\theta_nA(t)|}{v(t)}>\sup_{[\delta,1-\delta]}\frac{|u_n(F_n(t))+\sqrt n\theta_nA(t)|}{v(t)}\right)=o(1).
\eqno(A.15)
$$
Indeed, having (A.15), and arguing as in the proof of (ii) of Theorem 1, (5.7) imply that for positive $w$
$$P_{\theta_n}^n({\cal S}_n-\sqrt n\theta_nm_0\leqslant w)\geqslant { Pr}\left(\sup_{[\delta,1-\delta]}\frac{|u_n(F_n(t))|}{v(t)}\leqslant w\right)+o(1).$$
Since $u_n \circ F_n$ converges in distribution to a Brownian bridge, then (ii) follows. 

To prove (A.15) let $(\theta_n)$ be such that $n^{\varpi}\theta_n \to 0$ as $n \to \infty$. Let $(\iota_n)$ be a sequence such that $\iota_n\leq \log n$ and $\iota_n \to \infty$ as $n \to \infty$. Let $U_{(1)} \leq ... \leq U_{(n)}$ be order statistics of $n$ i.i.d. $U(0,1)$ random variables. Set
$$
\mathbb{E}_n = \Bigl\{\frac{1}{n \iota_n} \leq F_n^{-1}(U_{(1)}) \leq \frac{\iota_n}{n},\;1-\frac{\iota_n}{n} \leq F_n^{-1}(U_{(n)}) \leq 1 - \frac{1}{n \iota_n}\Bigr\}.
$$
Then, due to (5.6) and the assumption $n^{\varpi} \theta_n \to 0$, 
$$
\lim_{n \to \infty} Pr(\mathbb{E}_n) =1.
\eqno(A.16)
$$ 

Now, by the definitions of $m_0$ and $\delta$ in (5.7) and (A.16), we infer in the same way as in (A.13) that
$$ { Pr}\left(\sup_{(0,1)}\frac{|u_n(F_n(t))+\sqrt n\theta_nA(t)|}{v(t)}>\sup_{[\delta,1-\delta]}\frac{|u_n(F_n(t))+\sqrt n\theta_nA(t)|}{v(t)}\right)
$$
$$
\leqslant { Pr}\left(\sup_{(0,1)}\frac{|u_n(F_n(t))|}{v(t)}>\frac{m_0}{4}\sqrt n\theta_n\right)\leqslant 
{ Pr}\left(\Bigl\{\sup_{(0,1)}\frac{|u_n(F_n(t))|}{v(t)}>\frac{m_0}{4}\sqrt n\theta_n\Bigr\} \cap \mathbb{E}_n\right)+o(1).
\eqno(A.17)
$$
On the event $ \mathbb{E}_n$, for $t \in (0,F^{-1}_n (U_{(1)}))$ and $n$ sufficiently large, by (5.6) and $n^{\varpi} \theta_n \to 0$, it holds that 
$$
\frac{|u_n(F_n(t))|}{v(t)} \leq \sqrt n\sqrt{\frac{t}{1-t}} + \sqrt n \theta_n \frac{|A(t)|}{v(t)}\leq 2 \sqrt {\iota_n}.
\eqno(A.18)
$$
The same estimate holds on $ \mathbb{E}_n$ for $\;t \in (F^{-1}_n (U_{(n)}),1)$. On the other hand, on the event $\mathbb{E}_n $, for $t \in [F_n^{-1}(U_{(1)}),F_n^{-1}(U_{(n)})]$ and $n$ suficiently large, by (A.14) and (5.6), 
$$
\frac{|u_n(F_n(t))|}{v(t)} \leq \frac{|u_n(F_n(t))|}{v(F_n(t))}\sqrt{1 + \theta_n \frac{|A(t)|}{t(1-t)}} 
\leq 2 \frac{|u_n(F_n(t))|}{v(F_n(t))},
\eqno(A.19)
$$
provided that $\iota_n \to \infty$ is chosen in such a way that $(n \iota_n)^{\varpi} \theta_n \to 0$. 
The relations (A.18) and (A.19) allow to majorize the right hand side of (A.17) by
$$ 
Pr\Bigl(2 \sqrt{\iota_n} + 2 \sup_{(0,1)} \frac{|u_n(F_n(t)|}{v(F_n(t))} > \frac{m_0}{4}\sqrt n \theta_n \Bigr) + o(1)=Pr\Bigl(\sup_{(0,1)} \frac{|u_n(t)|}{v(t)} > \frac{m_0}{8}\sqrt n \theta_n -\sqrt{\iota_n}\Bigr) + o(1).
$$
In view of the assumption $(\log n \theta_n^2)/\log \log n \to \infty$ as $n \to \infty$, we have $\sqrt{\iota_n}/(\sqrt n\theta_n) \to 0$. An application of the main result of Mason (1985) concludes the proof of (ii). 

By (i) and (ii), ${\cal S}_n - \theta_n \sqrt n m_0$ is bounded in the probability $P_{\theta_n}$ and, in consequence, ${\cal S}_n /\theta_n \sqrt n m_0 \stackrel{P_{\theta_n}}{\longrightarrow} 1$. Hence, for $b_{\cal M}(P_{\theta_n}^n) = \sqrt { \log(\theta_n \sqrt n m_0)}$
it holds that ${\cal M}_n - b_{\cal M}(P_{\theta_n}^n) \stackrel{P_{\theta_n}}{\longrightarrow} 0$. 
Therefore, for $\{P_{\theta_n}\}$ satisfying the assumptions of Theorem 2, we infer that ${\cal M}_n /b_{\cal M}(P_{\theta_n}^n) \stackrel{P_{\theta_n}}{\longrightarrow} 1$.\hfill $\Box$\\

\noindent
{\bf A.6. Proof of Theorem 3}\\

\noindent
Let us start with an useful elementary result.\\

\noindent
{\bf Lemma A.1.} {\it Let $\{T_n\}$ be a sequence of non-negative random variables defined on a probability space with a measure $P$ which can depend on $n$. Moreover, let $\{\mu_n\}$ be a sequence of positive numbers tending to infinity as $n \to \infty$. 

\noindent
Then the following conditions are equivalent:

\vspace{2mm}
(i) $\;T_n - \mu_n \stackrel{ D}{\longrightarrow} T$;\\

(ii) $\displaystyle\frac{T_n^2 - \mu_n^2}{2\mu_n} \stackrel{D}{\longrightarrow} T$;\\
\noindent
and each of them implies $T_n/\mu_n \stackrel{P}{\longrightarrow} 1.$}\\

\noindent 
Set $T_n={\cal I}_n$ and $\mu_n=\sqrt{n}\theta_n||A^*||_2$. Recall that under $P_{\theta_n}^n$ the empirical process $\sqrt{n}(\hat{F}_n(t)-t)$ has the same distribution as $u_n(F_n(t))+\sqrt{n}\theta_nA(t)$, where $u_n(t)$ denotes the uniform empirical process. Hence
$$ 
\frac{T_n^2-\mu_n^2}{2\mu_n}=\frac{1}{2\sqrt{n}\theta_n||A^*||_2}d_n\left(\frac{u_n(F_n(t))}{[F_n(t)(1-F_n(t))]^{\ell}}\right)+\frac{1}{||A^*||_2}l_n\left(\frac{u_n(F_n(t))}{[F_n(t)(1-F_n(t))]^{\ell}}\right),
\eqno (A.20)
$$
where for $f\in D[0,1]$
$$
d_n(f)=\int_0^1 f^2(t)\left(1+\theta_n\frac{(1-2t)A(t)}{t(1-t)}-\theta_n^2\frac{A^2(t)}{t(1-t)}\right)^{2\ell} \frac{1}{[t(1-t)]^{1-2\ell}}dt,
$$
$$
l_n(f)=\int_0^1 f(t)\left(1+\theta_n\frac{(1-2t)A(t)}{t(1-t)}-\theta_n^2\frac{A^2(t)}{t(1-t)}\right)^{\ell} \frac{A(t)}{[t(1-t)]^{1-\ell}}dt.
$$
By the inequality $(a+b)^r\leq a^r+b^r, \;a,b>0,\;0<r<1,$ we have for sufficiently large $n$
$$
\left(1+\theta_n\frac{(1-2t)A(t)}{t(1-t)}-\theta_n^2\frac{A^2(t)}{t(1-t)}\right)^{r} \leq 1+\frac{|A(t)|^{r}}{[t(1-t)]^{r}}.
$$
Applying the above estimate to $r=2\ell$ and $r=\ell$, by (6.3) and the Lebesgue Dominated Theorem it follows
$$ 
d_n(f)\to d(f)=\int_0^1f^2(t)\frac{1}{[t(1-t)]^{1-2\ell}}dt,\;f\in D[0,1],
$$
$$
l_n(f)\to l(f)=\int_0^1 f(t)\frac{A(t)}{[t(1-t)]^{1-\ell}}dt,\;f\in D[0,1].
$$
Moreover, for any sequence $f_n(t)\in D[0,1]$ converging in $D[0,1]$ to $f\in C[0,1]$ it holds $d_n(f_n)\to d(f)$
and $l_n(f_n)\to l(f)$. As the process $B(t)/(t(1-t))^{\ell}$ has continuous trajectories a.s., then by Theorem 5.5 of Billingsley(1968) the above implies
$$
d_n\left(\frac{u_n(F_n(t))}{[F_n(t)(1-F_n(t))]^{\ell}}\right)\stackrel{\cal D}{\to}d\left(\frac{B(t)}{[t(1-t)]^{\ell}}\right)=\int_0^1 \frac{B^2(t)}{t(1-t)}dt
\eqno (A.21)
$$
and
$$
l_n\left(\frac{u_n(F_n(t))}{[F_n(t)(1-F_n(t))]^{\ell}}\right)\stackrel{\cal D}{\to}\int_0^1 B(t)\frac{A(t)}{t(1-t)}dt.
\eqno (A.22)
$$
The right hand side of (A.22) is a mean-zero Gaussian random variable with the variance $\rho^2_A$.
Using this, the assumption $n\theta_n^2\to\infty$ and (6.3) the proof follows from (A.20) - (A.22).\hfill $\Box$\\



\noindent
{\bf A.7. Proof of Theorem 5}\\

\noindent 
To prove (7.4) we shall exploit throughout the relation of ${\cal M}_n$ and ${\cal S}_n$, and corresponding results for ${\cal S}_n$. Take any $w >0$ and define $w_n^*=\sqrt {\log (1+ w + \theta_n \sqrt n m_0)}$. Set $\alpha_n = P_0^n({\cal M}_n \geq w_n^*)$. Since ${\cal M}_n$ has continuous and increasing distribution function then $t_{\alpha_n n} = w_n^*$ is the critical value of ${\cal M}_n$ corresponding to the level $\alpha_n$. 

By the assumption we have $(w_n^*)^2/\log \log n \geq [\log \theta_n\sqrt{n} + \log m_0]/ \log \log n\to\infty$ 
and Lemma 3 implies that
$$
-\frac{\log \alpha_n}{[w_n^*]^2} \to 2\;\;\;\mbox{as}\;\;n \to \infty.
$$
This yields $\alpha_n \to 0$ and $-[\log \alpha_n ] /n \to 0$ and means that $\{\alpha_n\}$ is an admissible significance level and the assumptions of Theorem 1 in Inglot et al. (2018) hold with $\gamma_n = \log \log n$ and $\lambda_n = n$.

On the other hand, the power of ${\cal M}_n$ under $P_{\theta_n}^n$ equals $P_{\theta_n}^n({\cal M}_n \geq w_n^*) = P_n^n({\cal S}_n - \sqrt n \theta_n m_0 \geq w)$. By Theorem 2 we have
$$
0<1-S_2(w) \leq \lim \inf_n P_{\theta_n}^n({\cal S}_n - \sqrt n\theta_n m_0 \geq w) 
\leq \lim \sup_n P_{\theta_n}^n({\cal S}_n - \sqrt n\theta_n m_0 \geq w) \leq 1-S_1(w) < 1.
$$
This proves that under $\{\alpha_n\}$ the test ${\cal M}_n$ has non-degenerate asymptotic power. By Proposition 1 of the present contribution and Theorem 1 in Inglot et al. (2018) the proof of (7.4) is concluded.\hfill $\Box$\\

\noindent
{\bf A.8. Proof of Lemma 7}\\

\noindent 
By (A.14) and the assumption $\theta_n^2m_n^2/\kappa_n\to 0$ we have for $t\in[\kappa_n,1-\kappa_n]$
$$0 \leq \frac{F_n(t)(1-F_n(t))}{t(1-t)}\leq 1+\frac{\theta_n\,m_n}{\sqrt{\kappa_n}}\frac{1-2\kappa_n}{\sqrt{1-\kappa_n}}= 1+o(1).
\eqno(A.23)$$

From the Chebyshev's inequality and (A.23) we get for arbitrary $\epsilon\in (0,1)$
$$P_{\theta_n}^n({\cal E}_n-\theta_n\sqrt{n}m_n\leq -\epsilon\theta_n\sqrt{n}m_n)\leq {Pr}\left(\frac{e_n(F_n(t_n))}{\sqrt{t_n(1-t_n)}}\leq -\epsilon\theta_n\sqrt{n}m_n\right)
\eqno(A.24)$$
$$\leq \frac{F_n(t_n)(1-F_n(t_n))}{t_n(1-t_n)}\frac{1}{\epsilon^2 n\theta_n^2m_n^2}= \frac{1}{\epsilon^2 n\theta_n^2m_n^2}(1+o(1))\to 0,$$
where $t_n$ is defined in (9.1). On the other hand, since $F_n(\kappa_n)\sim \kappa_n,\;F_n(1-\kappa_n)\sim 1-\kappa_n$, then by the triangle inequality for sufficiently large $n$ we have
$$P_{\theta_n}^n({\cal E}_n-\theta_n\sqrt{n}m_n\geq \epsilon\theta_n\sqrt{n}m_n)\leq {Pr}\left(\sup_{[\kappa_n,1-\kappa_n]}\frac{|e_n(F_n(t))|}{\sqrt{t(1-t)}}\geq \epsilon\theta_n\sqrt{n}m_n\right)
\eqno(A.25)$$
$$\leq {Pr}\left(\sup_{[\kappa_n/2,1-\kappa_n/2]}\frac{|e_n(t)|}{\sqrt{t(1-t)}}\geq \epsilon\theta_n\sqrt{n}m_n\right)$$ 
$$={Pr}\left(\sqrt{a_n}\sup_{[\kappa_n/2,1-\kappa_n/2]}\frac{|e_n(t)|}{\sqrt{t(1-t)}}-a_n-\frac{1}{2}\log (a_n/2\pi)\geq \sqrt{a_n}\epsilon\theta_n\sqrt{n}m_n-a_n-\frac{1}{2}\log (a_n/2\pi)\right),$$
where $a_n=2\log(\log(2/\kappa_n -1))$ is the normalizing sequence in the Darling-Erd\H{o}s theorem. The last expression tends to 0 due to the assumption $n\theta_n^2/\log\log(1/\kappa_n)\to\infty$ and the theorem of Jaeschke (1979).
Combining (A.24) and (A.25) we obtain (9.2).\hfill $\Box$
\\

\noindent
{\bf \large Acknowledgements}.} The work of B. \'Cmiel was partially supported by the Faculty of Applied Mathematics AGH UST dean grant for PhD students and young researchers within subsidy of Ministry of Science and Higher Education.\\


\begin{thebibliography}{99.}%

\bibitem{ } Adler, R.J. (1990). {\it An Introduction to Continuity, Extrema, and Related Topics for General Gaussian Processes}. Institute of Mathematical Statistics. Lecture Notes-Monograph Series. Vol. 12, Hayward, California.
\bibitem { } Anderson, T. W. and Darling, D. A. (1952). Asymptotic theory of certain ``goodness of fit'' criteria based on stochastic processes. {\it Ann. Math. Statist.} {\bf 23} 193-212.
\bibitem{ } Billingsley, P. (1968). {\it Convergence of Probability Measures}. Wiley.
\bibitem{ } Borovkov, A. A. and Sycheva, N. M. (1968). On asymptotically optimal non-parametric criteria. {\it Theory Probab. Appl.} {\bf 13} 359-393.
\bibitem{ } Borovkov, A. A. and Sycheva, N. M. (1970). On asymptotically optimal nonparametric criteria. In {\it Nonparametric Techniques in Statistical Inference} (M. L. Puri, ed.) 259-266. Cambridge Univ. Press.
\bibitem{ } Cai, T. T., J. Jeng and J. Jin (2011). Optimal detection of heterogeneous and heteroscedastic mixtures. {\it J. R. Statist. Soc. B} {\bf 73} 629-662.
\bibitem{ } Chicheportiche, R. and Bouchaud, J.-F. (2012). Weighted Kolmogorov-Smirnov test: Accounting for the tails. {\it Physical Review} {\bf E 86} 041115-1 - 041115-6.
\bibitem{ } Cont, R. (2001). Empirical properties of assets: stylized facts and statistical issues. {\it Quant. Fin.} {\bf 1} 223-236.
\bibitem{ } Chrampi, K. and Ycart, B. (2015). Weighted Kolmogorov Smirnov testing: and alternative for Gene Set Enrichment Analysis. {\it Stat. Appl. Genet. Mol. Biol.} {\bf 14} 279-293; {\it arXiv:1410.1620v1}.
\bibitem{ } Cs\"org\H{o}, M., Cs\"org\H{o}, S., Horvath, L. and Mason, D. M. (1986). Weighted empirical and quantile processes. {\it Ann. Probab.} {\bf 14}, 31-85.
\bibitem{ } Cs\"org\H{o}, M. and Horvath, L. (1993). {\it Weighted Approximations in Probability and Statistics}. Wiley.
\bibitem{ } Ditzhaus, M. (2018). Signal detection via Phi-divergences for general mixtures. {\it arXiv: 1803.06519v1}.
\bibitem {} Donoho, D. and Jin, J. (2004). Higher criticism for detecting sparse heterogeneous mixtures. {\it Ann. Statist.} {\bf 32} 962-994.
\bibitem {} Donoho, D. and Jin, J. (2015). Higher criticism for large-scale inference, especially for rare and weak effects. {\it Stat. Sci.} {\bf 30}, 1-25.
\bibitem{ } Eicker, F. (1979). The asymptotic distribution of the suprema of the standardized empirical process. {\it Ann. Statist.} {\bf 7} 116-138.
\bibitem{ } Ermakov, M.S. (2004). On asymptotically efficient statistical inference for moderate deviation probabilities. {\it Theory Probab. Appl.} {\bf 48}, 622-641.
\bibitem { } Gontscharuk, V., Landwehr, S. and Finner, H. (2016). Goodness of fit tests in terms of local levels with special emphasis on higher criticism tests. {\it Bernoulli} {\bf 22} 1331-1363.
\bibitem { } Grabchak, M. and Samorodnitsky, G. (2010). Do financial returns have finite or infinite variance? A paradox and an explanation. {\it Quant. Finance} {\bf 10} 883-893.
\bibitem { } Greenshtein, E. and Park, J. (2012). Robust test for detecting a signal in a high dimensional sparse normal vector. {\it J. Statist. Plann. Inference} {\bf 142} 1445-1456.
\bibitem { } Handcock, M. S. and Morris, M. (1999). {\it Relative Distribution Methods in the Social Sciences.} Springer, New York.
\bibitem{ } Inglot, T. (1999). Generalized intermediate efficiency of goodness of fit tests. {\it Math. Methods Statist.} {\bf 8} 487-509.
\bibitem{ } Inglot, T. (2019). Intermediate efficiency of tests under heavy-tailed alternatives. {\it arXiv:1902.06622v1 [math.ST]}.
\bibitem{ } Inglot, T., Kallenberg, W.C.M. and Ledwina, T. (1998). Vanishing shortcoming and asymptotic relative efficiency. {\it Memorandum} 1467, Faculty of Mathematical Sciences, Univ. Twente.
\bibitem{ } Inglot, T., Kallenberg, W.C.M. and Ledwina, T. (2000). Vanishing shortcoming and asymptotic relative efficiency. {\it Ann. Statist.} 28, 215-238.
\bibitem { } Inglot, T. and Ledwina, T. (1990). On probabilities of excessive deviations for Kolmogorov-Smirnov, Cram{\'e}r-von Mises and chi-square statistics. {\it Ann. Statist.} {\bf 18} 1491-1495. 
\bibitem { } Inglot, T. and Ledwina, T. (1993). Moderately large deviations and expansions of large deviations for some functionals of weighted empirical process. {\it Ann. Probab}. {\bf 21} 1691-1705.
\bibitem{ } Inglot, T. and Ledwina, T. (1996). Asymptotic optimality of data driven Neyman's tests for uniformity. {\it Ann. Statist.} {\bf 24} 1982-2019.
\bibitem{ } Inglot, T. and Ledwina, T. (2006). Intermediate efficiency of some max-type statistics. {\it J. Statist. Plan. Inference.} {\bf 136} 2918-2935.
\bibitem{ } Inglot, T., Ledwina, T. and \'Cmiel, B. (2018). Intermediate efficiency in nonparametric testing problems with an application to some weighted statistics. {\it ESAIM: Probab. Stat.} accepted; for a preliminary variant see {\it arXiv:1806.02020v1 [math.ST]}.
\bibitem{ } Ingster, Y. (1997). Some problems of hypothesis testing leading to infinitely divisible distributions. {\it Math. Methods Statist.} {\bf 6}, 47-69.
\bibitem { } Ivchenko, G. I. and Mirakhmedov Sh. A. (1995). Large deviations and intermediate efficiency of decomposable statistics in a multinomial scheme. {\it Math. Methods Statist.} {\bf 4} 294-311.
\bibitem { } Jager, L. and Wellner, J. A. (2004). On the ``Poisson boundaries'' of the family of weighted Kolmogorov statistics. In \textit{Festschrift for Herman Rubin} (A. DasGupta, ed.) 319-331. IMS, Beachwood, OH.
\bibitem { } Jager, L. and Wellner, J. A. (2007). Goodness-of-fit tests via phi-divergences. {\it Ann. Statist.} {\bf 35} 2018-2053.
\bibitem { } Janssen A. (1995). Principal component decomposition of non-parametric tests. {\it Probab. Theory Relat. Fields} {\bf 101} 193-209.
\bibitem { } Jin, J., Starck, J. -L., Donoho, D. L., Aghanim, N., and Forni, O. (2005), Cosmological non-Gaussian signature detection: Comparing performance of different statistical tests, {\it EURASIP J. Appl. Signal Processing} {\bf 15} 2470-2485.
\bibitem{ } Kallenberg, W. C. M. (1983). Intermediate efficiency, theory and examples. {\it Ann. Statist.} {\bf 11} 1401-1420.
\bibitem{ } Lehmann, E. L. (1953). The power of rank tests. {\it Ann. Math. Staist.} {\bf 24} 23-43.
\bibitem{ } Ledwina, T. and Wy{\l}upek, G. (2012). Nonparametric tests for first order stochastic dominance. {\it TEST} {\bf 21}, 730-756.
\bibitem{ } Li, J. and Siegmund D. (2015). Higher criticism: $p$-values and criticism. {\it Ann. Statist.} {\bf 43} 1323-1350.
\bibitem{ } Lifshits, M. A. (1995). {\it Gaussian Random Functions}. Springer Science+Business Media, Dordrecht.
\bibitem { } Lockhart, R. A. (1991). Overweight tails are inefficient. {\it Ann. Statist.} {\bf 19} 2254-2258.
\bibitem { } Mason, D. M. (1985). Some large deviation results for weighted empirical processes. {\it Statistics $\&$ Decisions} {\bf 2} 89-98.
\bibitem{ } Mason, D. M. and Schuenemeyer, J. H. (1983). A modified Kolmogorov-Smirnov test sensitive to tail alternatives. {\it Ann. Statist.} {\bf 11} 933-946.
\bibitem { } Mason, D. M. and Schuenemeyer, J. H. (1992). Correction: A modified Kolmogorov-Smirnov test sensitive to tail alternatives. {\it Ann. Statist.} {\bf 20} 620-621.
\bibitem { } Merlev\`ede, F. and Peligrad, M. (2009). Functional moderate deviations for triangular arrays and applications. {\it ALEA Lat. Am. J. Probab. Math. Stat.} {\bf 5} 3-20.
\bibitem { } Milbrodt, H. and Strasser, H. (1990). On the asymptotic power of the two-sided Kolmogorov-Smirnov test. {\it J. Statist. Plann. Inference } {\bf 26} 1-23.
\bibitem{ } Moscovich, A., Nadler, B. and Spiegelman, C. (2016). On the exact Berk-Jones statistics and their $p$-value calculation. {\it Electron. J. Stat.} {\bf 10} 2329-2354.
\bibitem{ } Pearson, E. S., D'Agostino, R. B. and Bowman, K. O. (1977). Tests for departures from normality: Comparison of powers. {\it Biometrika} {\bf 64} 231-246.
\bibitem{ } R\'ev\'esz, P. (1982). A joint study of the Kolmogorov-Smirnov and the Eicker-Jaeschke statistics. {\it Statistics $\&$ Decisions} {\bf 1}, 57-65.
\bibitem { } Stepanova, N. and Pavlenko, T. (2018). Goodness-of-fit tests based on sup-functionals of weighted empirical processess. {\it Theory Probab. Appl.} {\bf 63}, 292-317.
\bibitem { } Thas, O. (2010). {\it Comparing Distributions}. Springer, New York.
\bibitem{ } van der Vaart, W. and Wellner, J.A. (2000). {\it Weak Convergence and Empirical Processes}. Springer Series in Statistics.
 
\end{thebibliography}
\end{document}